\documentclass[10pt,twoside]{article}
\usepackage[margin=1in]{geometry}
\usepackage{latexsym,amssymb,amsmath,amsfonts,xcolor,cases,mathtools}
\usepackage{graphicx,theorem}
\usepackage[caption=false]{subfig}
\usepackage{chngcntr}
\usepackage{natbib}
\usepackage{epstopdf}
\usepackage{hyperref}
\usepackage{lipsum}
\usepackage{algorithmic}
\usepackage{amsopn}
\usepackage{fancyhdr}

\numberwithin{equation}{section}

\hypersetup{
    colorlinks = true,
    linkcolor = {blue},
    citecolor = {blue},
    urlcolor  = {blue},}

\newcommand{\RE}[1]{\textup{Re}\left\{#1\right\}}
\newcommand{\IM}[1]{\textup{Im}\left\{#1\right\}}
\newcommand{\CC}{{\mathbb{C}}}

\newcommand{\ZZ}{{\mathbb{Z}}}

\newcommand{\Der}[2]{\frac{\text{d}#1}{\text{d}#2}}

\newcommand{\BF}[1]{\mathbf{#1}}

\newcommand{\RED}{} 

\begin{document}
\title{Diffraction of acoustic waves by a wedge of point scatterers\footnotemark[1]}
\author{M. A. Nethercote\footnotemark[2],~A. V. Kisil\footnotemark[2]~and~R. C. Assier\footnotemark[2]}
\footnotetext[1]{Submitted to the SIAM Journal on Applied Mathematics on 06/08/2021. \newline 
Funding: A.V.K. is supported by a Royal Society Dorothy Hodgkin Research Fellowship and a Dame Kathleen Ollerenshaw Fellowship, through which M.A.N. has been funded.}
\footnotetext[2]{Department of Mathematics, Alan Turing Building, University of Manchester, United Kingdom,\newline
		(matthew.nethercote@manchester.ac.uk, anastasia.kisil@manchester.ac.uk, raphael.assier@manchester.ac.uk).}

\date{}
\maketitle

\begin{abstract}
This article considers the problem of diffraction by a wedge consisting of two semi-infinite periodic arrays of point scatterers. The solution is obtained in terms of two coupled systems, each of which is solved using the discrete Wiener--Hopf technique. An effective and accurate iterative numerical procedure is developed to solve the diffraction problem, which allows us to compute the interaction of thousands of scatterers forming the wedge. A crucial aspect of this numerical procedure is the effective truncation of slowly convergent single and double infinite series, which requires careful asymptotic analysis. A convergence criteria is formulated and shown to be satisfied for a large class of physically interesting cases. A comparison to direct numerical simulations is made, highlighting the accuracy of the method.

%
\end{abstract}

%

\section{Introduction}\label{Intro}
Wedge diffraction is a classical and well-studied problem \cite{WedgeReview}. It has been considered in different media and with many different boundary conditions (for hard, soft, impedance, penetrable and elastic wedges). Diffraction by periodic structures is also an interesting and challenging problem, which has been studied mainly for infinite~\cite{Bruno_Fernandez17, MovchanBook, Scarpetta_Sumbatyan96} and semi-infinite arrays \cite{HillsKarp1965,LintonMartin2004}. Much less is known about diffraction by wedges made from periodic structures (in particular point scatterers), which is the object of this article.

There has been renewed interest in periodic structures due to the growth of the metamaterials research. Frequently in metamaterial research, some form of periodicity is assumed which allows us to restrict the study to a small portion known as the unit cell. The global behaviour is then reconstructed from this local behaviour by symmetry. By modifying the unit cell, it is possible to modify the properties of the whole material drastically. It can display abnormal phenomena such as negative refraction, superlensing and backward wave propagation which allow for the control and manipulation of waves in ways that were originally perceived as unphysical. This has led to many studies, which increase the complexity of the unit cell and reconstruct the global scattering using the periodicity of the metamaterial~\cite{MovchanBook,Craster-book}.

An alternative approach looks into the case where periodicity is no longer assumed. This means that the metamaterial will have well-defined boundaries that can symbolise many different interfaces such as edges and corners. In this framework, this article could be considered as diffraction by the simplest possible wedge made out of a metasurface. Each face will be a semi-infinite array made of sound-soft cylinders with infinite height and a small radius, see Figure~\ref{fig:HW-diagram}.

The wedge configuration can also be viewed as two separate semi-infinite arrays with two sets of scattering coefficients to determine. From the two derived infinite system of equations, we will find a solution using the discrete Wiener--Hopf technique. The key steps of the technique are the application of the Z-transform and the factorisation of the Wiener--Hopf kernel (an infinite sum with some resemblance to the periodic Green's function) on the unit circle~\cite{CHCA2019}. Numerically, a number of challenges need to be overcome due to slow converging series. Some tools that are employed are tail-end asymptotics and rational approximations. It is found that the latter approach captures the essential interesting behaviour and is computationally fast. An iterative scheme is then constructed from the individual Wiener--Hopf solutions.

One interesting feature of periodic structures is the possibility of resonance effects. In the context of array scattering, it is frequently referred to as Wood's anomaly~\cite{Liu_Declercq15}. It was first observed in diffraction gratings experimentally \cite{Wood1902} and partly discussed theoretically by Lord Rayleigh \cite{Rayleigh1907}. At resonance, the solution usually takes a simple form of a plane wave moving along the array. Also, typically around the resonant state, interesting phenomena can be observed such as, for example, abnormal transmission or reflection~\cite{Korolkov_Nazarov16, Shanin_Korolkov17}. There are also different notions of resonance associated to arrays (and no consistent terminology), some of which are the resonances of an infinite array \cite{Evans_Porter_99} like Wood’s anomaly and some are associated with effects of finiteness of the array~\cite{Bennetts_Malte_17}. A wedge with periodic perturbations was considered in \cite{Kamotskii_Nazarov99_I} and \cite{Nazarov08}. In this article, we will not  be considering these cases since they have distinctive features which require special attention and will be considered in our future work.

The structure of the article is as follows. In Section \ref{IG}, we start by reviewing some fundamental canonical scattering problems related to infinite and semi-infinite arrays. In particular, we review several important aspects of multiple scattering theory, including Foldy's approximation, isotropic scattering and the discrete Wiener--Hopf technique. Section \ref{HW} focuses entirely on diffraction by a wedge of point scatterers. After formulating the problem, we derive an implicit solution via the discrete Wiener--Hopf technique and find an approximation using iterative methods. This is followed by discussions on the conditional convergence of the iterative scheme and comparisons with the finite element software COMSOL MultiPhysics version 5.5. Lastly, in Section \ref{ConFW} we draw some conclusions and discuss further work. The appendix describes the numerical methods, tail-end asymptotics and rational approximations used to make the proposed methods fast and reliable computationally.

\section{Scattering by infinite and semi-infinite arrays}\label{IG}
\begin{figure}[ht]\centering
\includegraphics[width=0.9\textwidth]{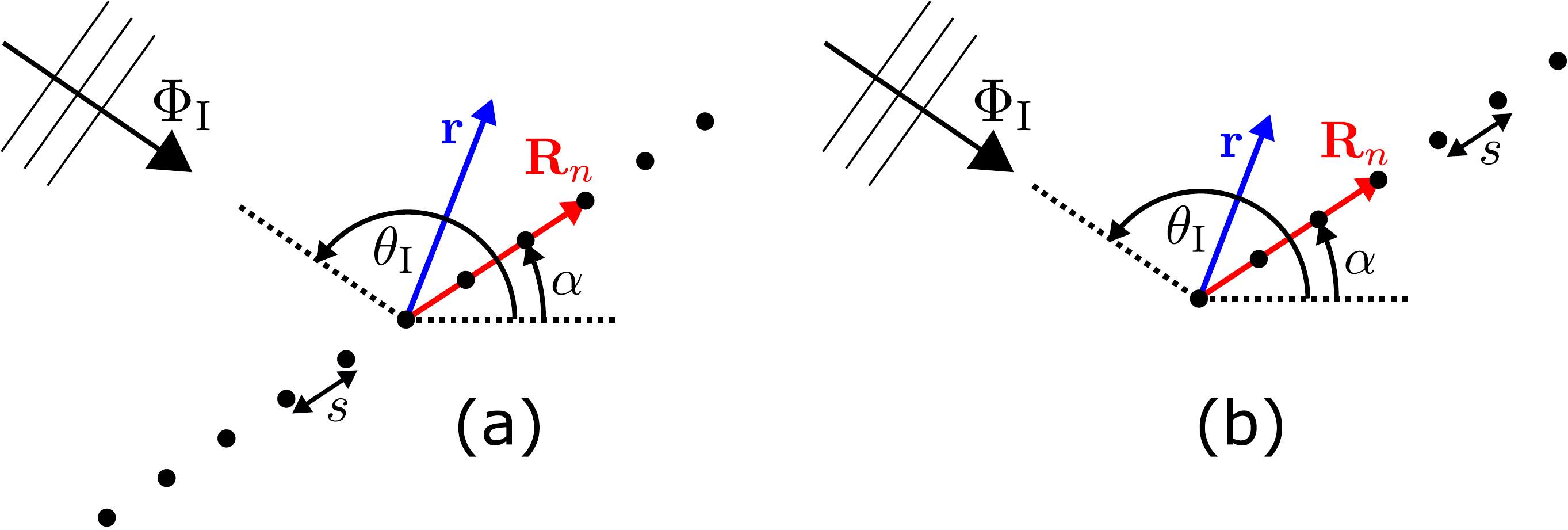}
\caption{Diagrams of the infinite (a) and semi-infinite (b) arrays with the cylinders located at $\BF{R}_n$, the position vector $\BF{r}$ and the incident wave $\Phi_{\text{I}}$.}
\label{fig:inf-array-diagram}
\end{figure}

In this article, we are investigating multiple scattering problems involving semi-infinite arrays comprised of equidistant scatterers. The positions of the scatterer centres are denoted by $\BF{R}_n(s,\alpha)$ where $n\in\ZZ$ is the scatterer index, $s>0$ is the distance between scatterers and $\alpha\in[0,\pi)$ is the angle between the array and the $x$-axis. The explicit definition of $\BF{R}_n$ will change throughout this article (see Figure \ref{fig:inf-array-diagram} and \ref{fig:HW-diagram}). The scatterers themselves are all cylinders of infinite height and small radius $a\in[0,s/2]$ and satisfy homogeneous Dirichlet boundary conditions. The upper limit of $a$ is to prevent overlapping cylinders. The incident wave field takes the form of a unit amplitude plane wave at normal (non-skew) incidence. Problems of this type are invariant with respect to the height position on the cylinder scatterers, which simplifies the three-dimensional problem to a two-dimensional projection.

We shall also be looking for time-harmonic solutions to the linear wave equation by assuming and then suppressing the time factor $e^{-i\omega t}$ (where $\omega$ is the angular frequency). Let $\Phi$ be the total wave field, which is decomposed into an incident field $\Phi_{\text{I}}$ and the resulting scattered field $\Phi_{\text{S}}$ by the equation $\Phi=\Phi_{\text{I}}+\Phi_{\text{S}}$ and use a polar coordinate system $(r,\theta)$ with the position vector $\BF{r}$. The incident wave field $\Phi_{\text{I}}(\BF{r})$ is given by,
\begin{align}\label{IG-incwave}
\Phi_{\text{I}}(\BF{r})=e^{-ikr\cos(\theta-\theta_{\text{I}})},
\end{align}
where $k$ is the wavenumber and $\theta_{\text{I}}$ is the incoming incident angle. 

We assume that the cylinder scatterers are widely spaced out compared to their size i.e. $a\ll s/2$ and we also assume that the cylinder radius is small in comparison to the wavelength of the incident wave ($a\ll 2\pi/k$). These assumptions alongside the Dirichlet boundary conditions allow us to make the assumption of isotropic scattering in the spirit of \cite{Foldy1945}. This means that for every cylinder, the scattered wave can be approximated by a Hankel function of the first kind and zeroth order (denoted $H^{(1)}_0$) multiplied by an unknown scattering coefficient. The total scattered field $\Phi_{\text{S}}$ (which satisfies Helmholtz's equation) is therefore a summation of all these individual scattered waves,
\begin{align}\label{IG-gensol-phis}
\Phi_{\text{S}}(\BF{r})=\sum_{n} \left[A_nH^{(1)}_0(k|\BF{r}-\BF{R}_n|)\right],
\end{align}
where $A_n$ are the unknown scattering coefficients and the summation limits will depend on the specifics of the array layout.

We can define the wave field incident on the $m^{\text{th}}$ scatterer, $\Phi_{\text{I}}^{(m)}$, as the original incident wave plus the sum of all the scattered waves from the rest of the cylinders. In other words, it is the total wave field minus the $m^{\text{th}}$ scattered wave,
\begin{align}\label{IG-incwave_m}
\nonumber\Phi_{\text{I}}^{(m)}(\BF{r})&=\Phi_{\text{I}}(\BF{r})+\sum_{n\neq m} \left[A_nH^{(1)}_0(k|\BF{r}-\BF{R}_n|)\right],\\
&=\Phi(\BF{r})-A_mH^{(1)}_0(k|\BF{r}-\BF{R}_m|).
\end{align}
After applying the Dirichlet boundary conditions to \eqref{IG-incwave_m}, we use the assumption that the cylinder radius $a$ is a small parameter to approximate, 
\begin{align}\label{IG-Am-BCs}
A_m=-\frac{\Phi_{\text{I}}^{(m)}(\BF{R}_m+a\BF{\hat{r}})}{H^{(1)}_0(ka)}\approx-\frac{\Phi_{\text{I}}^{(m)}(\BF{R}_m)}{H^{(1)}_0(ka)}.
\end{align}
Note that normally $\Phi$ cannot be evaluated at $\BF{R}_m$ because it diverges, but $\Phi_{\text{I}}^{(m)}$ can since the singularity has been removed. We can use this to construct an infinite system of equations to solve for the scattering coefficients,
\begin{align}\label{IG-general-system}
A_mH^{(1)}_0(ka)+\sum_{n\neq m} \left[A_nH^{(1)}_0(k|\BF{R}_m-\BF{R}_n|)\right]+\Phi_{\text{I}}(\BF{R}_m)=0.
\end{align} 

Before investigating the wedge problem in Section \ref{HW}, we will review two well-studied problems: the infinite array (see Figure \ref{fig:inf-array-diagram}(a)) and the semi-infinite array (see Figure \ref{fig:inf-array-diagram}(b)). The scatterers of both problems have locations given by,
\begin{align}\label{IG-Rn}
\BF{R}_n=ns(\cos(\alpha)\BF{\hat{x}}+\sin(\alpha)\BF{\hat{y}}),\ \text{where}\ 
\begin{cases}
n=0,\pm1,\pm2,...&\text{infinite array,}\\
n=0,1,2, ...&\text{semi-infinite array,}
\end{cases}
\end{align}
in the \RED{C}artesian basis $(\BF{\hat{x}},\BF{\hat{y}})$.

For the infinite array problem (denoted with superscript `inf'), we have the following infinite system of equations,
\begin{align}\label{IG-Am-system}
A^{\text{(inf)}}_mH^{(1)}_0\!(ka)+\!\!\sum_{\substack{n=-\infty \\n\neq m}}^\infty\! \left[A^{\text{(inf)}}_nH^{(1)}_0\!(ks|m-n|)\right]\!+e^{-iksm\cos(\theta_{\text{I}}-\alpha)}\!=\!0,\ \ m\in\ZZ.
\end{align} 
Since $\Phi_{\text{I}}(\BF{r}+\BF{R}_m)=e^{-iksm\cos(\theta_{\text{I}}-\alpha)}\Phi_{\text{I}}(\BF{r})$, it is logical to say that $\Phi$ (and hence $A^{\text{(inf)}}_m$) satisfies the same type of quasi-periodicity. We can therefore use $A^{\text{(inf)}}_m=e^{-iksm\cos(\theta_{\text{I}}-\alpha)}A^{\text{(inf)}}_0$ so that \eqref{IG-Am-system} is reduced to solving for a single coefficient $A^{\text{(inf)}}_0$ which is then given by,\footnote{\RED{The formula \eqref{IG-A0-SchSeries} is very slow to converge, which is why, in practice, we sought an alternative formula by using the same identities that lead to \eqref{Kernel-K(z)-fast}. The result is also quite similar to equation (2.47) in \cite{Linton1998}}}
\begin{align}\label{IG-A0-SchSeries}
A^{\text{(inf)}}_0=-\left(H^{(1)}_0(ka)+2\sum_{n=1}^\infty \left[\cos(ksn\cos(\theta_{\text{I}}-\alpha))H^{(1)}_0(ksn)\right]\right)^{-1}.
\end{align}

Similarly, for the semi-infinite array problem (with superscript `s-inf'), we have $A^{\text{(s-inf)}}_m$ satisfying the following infinite system of equations,
\begin{align}\label{SIG-Am-system}
A^{\text{(s-inf)}}_mH^{(1)}_0\!(ka)+\!\!\sum_{\substack{n=0 \\n\neq m}}^\infty\! \left[A^{\text{(s-inf)}}_nH^{(1)}_0\!(ks|m-n|)\right]\!+e^{-iksm\cos(\theta_{\text{I}}-\alpha)}\!=\!0,\ \ m\geq 0,
\end{align} 
and $A^{\text{(s-inf)}}_m=0$ for all $m<0$. 

We solve this system using the discrete analogue of the Wiener--Hopf technique. Here, we will only give an overview of the technique for brevity (see \cite{HillsKarp1965} and \cite{LintonMartin2004} for more details). We start by extending the system of equations \eqref{SIG-Am-system} for negative $m$ using some unknown coefficients $F_m$,
\begin{align}\label{SIG-Am-before-trans}
A^{\text{(s-inf)}}_mH^{(1)}_0\!(ka)+\!\!\sum_{\substack{n=-\infty\\n\neq m}}^\infty\! \left[A^{\text{(s-inf)}}_nH^{(1)}_0\!(ks|m-n|)\right]\!=\!
\begin{cases}-e^{-iksm\cos(\theta_{\text{I}})}&m\geq0,\\F_m&m<0.\end{cases}
\end{align}
Then we multiply \eqref{SIG-Am-before-trans} by $z^m$ and summing from $m=-\infty$ to $\infty$ (i.e. the discrete $Z$-transform) to obtain the Wiener--Hopf equation, 
\begin{align}\label{SIG-WHE}
K(z)A^+(z)=F^+(z)+F^-(z),
\end{align}
where the goal is to determine the $Z$-transform of the scattering coefficients $A^+(z)$,
\begin{align}\label{SIG-Z-transform}
A^+(z)=\sum_{m=-\infty}^\infty A^{\text{(s-inf)}}_mz^m=\sum_{m=0}^\infty A^{\text{(s-inf)}}_mz^m.
\end{align}
To help with the convergence of the $Z$-transform, it is useful to assume that $k$ has a small positive imaginary part (i.e. $0<\IM{k}\ll1$). Note that we will let $\IM{k}\rightarrow0$ after we finish using the Wiener--Hopf technique because the final solution is still valid when $\IM{k}=0$. This will be consistent with all future numerical computations presented in this article. 

We use the superscripts `$+$' (resp. `$-$') for functions that are analytic and zero-free in the region $\Omega^+$ (resp. $\Omega^-$), where, 
\begin{align}\label{SIG-regions}
\Omega^+=\left\{z\in\CC: |z| < e^{-\IM{k}s\cos(\theta_{\textrm{I}}-\alpha)}\right\},\ \ 
\Omega^-=\left\{z\in\CC: |z| > e^{-\IM{k}s}\right\}.
\end{align}
Note that the union $\Omega^+\cup\Omega^-$ is the entire complex $z$-plane and the intersection $\Omega^+\cap\Omega^-$ is the annulus region $e^{-\IM{k}s}<|z|< e^{-\IM{k}s\cos(\theta_{\textrm{I}}-\alpha)}$. In the Wiener--Hopf equation \eqref{SIG-WHE}, $F^+(z)$ is the discrete $Z$-transform of the incident wave forcing and $F^-(z)$ is an unknown forcing term:
\begin{align}\label{SIG-F+}
F^+(z)=\frac{1}{ze^{-iks\cos(\theta_{\text{I}}-\alpha)}-1},\quad F^-(z)=\sum_{m=-\infty}^{-1} F_mz^m,
\end{align}
where the simple pole at $z=e^{iks\cos(\theta_{\text{I}}-\alpha)}$ is considered to be outside the region $\Omega^+$. $K(z)$ is the discrete Wiener--Hopf kernel given by,
\begin{align}\label{SIG-K(z)}
K(z)=H^{(1)}_0(ka)+\sum_{\ell=1}^\infty \left[(z^\ell+z^{-\ell})H^{(1)}_0(ks\ell)\right].
\end{align}
The Wiener--Hopf kernel has the property $K(z)=K(1/z)$ for all $z$ and is both analytic and zero-free on the annulus region $e^{-\IM{k}s}<|z|< e^{\IM{k}s}$ (which includes the intersection region $\Omega^+\cap\Omega^-$) except for the branch point singularities at $z=e^{\pm iks}$. Note that \eqref{SIG-K(z)} is not well-suited for numerical evaluation. This is because the formula only converges on $e^{-\IM{k}s}<|z|< e^{\IM{k}s}$ and converges very slowly (see \eqref{Kernel-K-term-asymp}). In Appendix \ref{App-Kernel}, we will discuss some alternative methods to evaluate $K(z)$ including a different formula \eqref{Kernel-K(z)-fast} and an approximation \eqref{Kernel-K_approx-basic} which will be essential for producing the results of this article.

To proceed with the Wiener--Hopf technique, \eqref{SIG-WHE} needs to be rearranged such that each side of the equality is both analytic and zero-free in $\Omega^+$ or $\Omega^-$. A key step of the technique is the factorisation of the kernel, $K(z)=K^+(z)K^-(z)$, which is done such that the factors satisfy the property $K^+(z)=K^-(1/z)$ for all $z$. This factorisation is achieved by using Cauchy's integral formula to sum split the natural logarithm of the kernel
\begin{align}\label{SIG-CIs}
\ln(K(z))&=\frac{1}{2\pi i}\int_{C^+}\frac{\ln(K(\xi))}{\xi-z}\text{d}\xi-\frac{1}{2\pi i}\int_{C^-}\frac{\ln(K(\xi))}{\xi-z}\text{d}\xi,
\end{align}
where the integration contours $C^\pm$ are anticlockwise circular paths contained inside $\Omega^+\cap\Omega^-$ on the $\xi$ complex plane. Additionally $C^+$ (resp. $C^-$) will also run radially above (resp. below) the pole at $\xi=z$. The former (resp. latter) integral is analytic in $\Omega^+$ (resp. $\Omega^-$). If we apply a simple transformation ($\xi\rightarrowtail 1/\xi$) to the $C^+$ integral, we can rewrite \eqref{SIG-CIs} such that the factors will obviously satisfy $K^+(z)=K^-(1/z)$,
\begin{align}
\label{SIG-K^+(z)}\ln(K^+(z))&=\ln(K_0)-\frac{1}{2\pi i}\int_{C^-}\frac{z\ln(K(\xi))}{z\xi-1}\text{d}\xi,\\
\label{SIG-K^-(z)}\ln(K^-(z))&=\ln(K_0)-\frac{1}{2\pi i}\int_{C^-}\frac{\ln(K(\xi))}{\xi-z}\text{d}\xi,\\
\label{SIG-K0}\text{where,}\quad \ln(K_0)&=\ln(K^+(0))=\frac{1}{4\pi i}\int_{C^-}\frac{\ln(K(\xi))}{\xi}\text{d}\xi.
\end{align}
We use pole removal to sum split $\frac{F^+(z)}{K^-(z)}$ after factorising the kernel and then apply Liouville's theorem to solve the Wiener--Hopf equation to obtain the transformed solution, 
\begin{align}\label{SIG-A(z)}
A^+(z)=\frac{1}{K^-(e^{iks\cos(\theta_{\text{I}}-\alpha)})K^+(z)(ze^{-iks\cos(\theta_{\text{I}}-\alpha)}-1)}.
\end{align}

Now, we apply the identity $K^-(z)=K^+(1/z)$ and the inverse $Z$-transform,
\begin{align}\label{SIG-Z-inverse}
A^{\text{(s-inf)}}_m=\frac{1}{2\pi i}\int_{C^-} A^+(z)z^{-m-1}\text{d}z,
\end{align}
to determine the scattering coefficients for the semi-infinite array problem,
\begin{align}\label{SIG-Am-final}
A^{\text{(s-inf)}}_m=-\frac{1}{K^+(e^{-iks\cos(\theta_{\text{I}}-\alpha)})}\sum_{n=0}^m\left[\lambda_ne^{-iks(m-n)\cos(\theta_{\text{I}}-\alpha)}\right].
\end{align}
Here the coefficients $\lambda_n$ are defined by the following Taylor series,
\begin{align}\label{SIG-1/K-series}
\frac{1}{K^+(z)}=\sum_{n=0}^\infty \lambda_n z^n,\ \ \text{where}\ \ \lambda_n=\frac{1}{n!}\Der{^n}{z^n}\left[\frac{1}{K^+(z)}\right]_{z=0},
\end{align}
and can be given using Cauchy's integral formula,
\begin{align}\label{SIG-lambda-def}
\lambda_n=\frac{1}{2\pi i}\int_{C^-}\frac{z^{-n-1}}{K^+(z)}\text{d}z=\frac{1}{2\pi}\int_{-ks}^{2\pi-ks}\frac{e^{-int}}{K^+(e^{it})}\text{d}t.
\end{align}
It should be noted that $\frac{1}{K^+(e^{it})}$ is $L^1$ integrable on the real line segment $[-ks,2\pi-ks]$ (i.e. the integral $\int_0^{2\pi}\left|\frac{1}{K^+(e^{it})}\right|\text{d}t$ is finite). By the Riemann-Lebesgue lemma, we can therefore say that $\lambda_n\rightarrow0$ as $n\rightarrow\infty$. The $\lambda_n$ coefficients for $n\geq1$ can also be written without requiring an exact formula for the $K^+$ factor. We do this by using the original definition \eqref{SIG-1/K-series}, the following identity,
\begin{align}\label{SIG-1st-der-con}
\Der{}{z}\left[\frac{1}{K^+(z)}\right]=-\frac{1}{K^+(z)}\Der{}{z}\left[\ln\left(K^+(z)\right)\right],
\end{align}
and the derivatives of \eqref{SIG-K^+(z)}, which results in the following formula, 
\begin{align}\label{SIG-lambda_n-integral}
\lambda_n=-\sum_{m=1}^n\frac{m\lambda_{n-m}}{n}\int_{0}^1 \cos(m\pi\tau)\ln\left(K(e^{i\pi\tau})\right)\text{d}\tau,\quad n\geq1.
\end{align}

Lastly, we can also rewrite the equation for the scattering coeffcients \eqref{SIG-Am-final} as a recurrence relation,
\begin{align}\label{SIG-Am-recurrence}
A^{\text{(s-inf)}}_m&\!=\!e^{-iks\cos(\theta_{\text{I}}-\alpha)}A^{\text{(s-inf)}}_{m-1}\!+\!\frac{\lambda_m}{\lambda_0}A^{\text{(s-inf)}}_0,\ \  A^{\text{(s-inf)}}_0\!=\!-\frac{\lambda_0}{K^+(e^{-iks\cos(\theta_{\text{I}}-\alpha)})}\cdot
\end{align}
This recurrence relation is useful because it is more efficient for computing the scattering coefficients $A^{\text{(s-inf)}}_m$ and it highlights their asymptotic behaviour as $m\rightarrow\infty$. Specifically, they become equivalent to the solution of the infinite array problem $A^{\text{(s-inf)}}_m\sim A^{\text{(inf)}}_m$ as $m\rightarrow\infty$.

\section{Diffraction by a wedge of point scatterers}\label{HW}
\begin{figure}[h!]\centering
\includegraphics[width=0.5\textwidth]{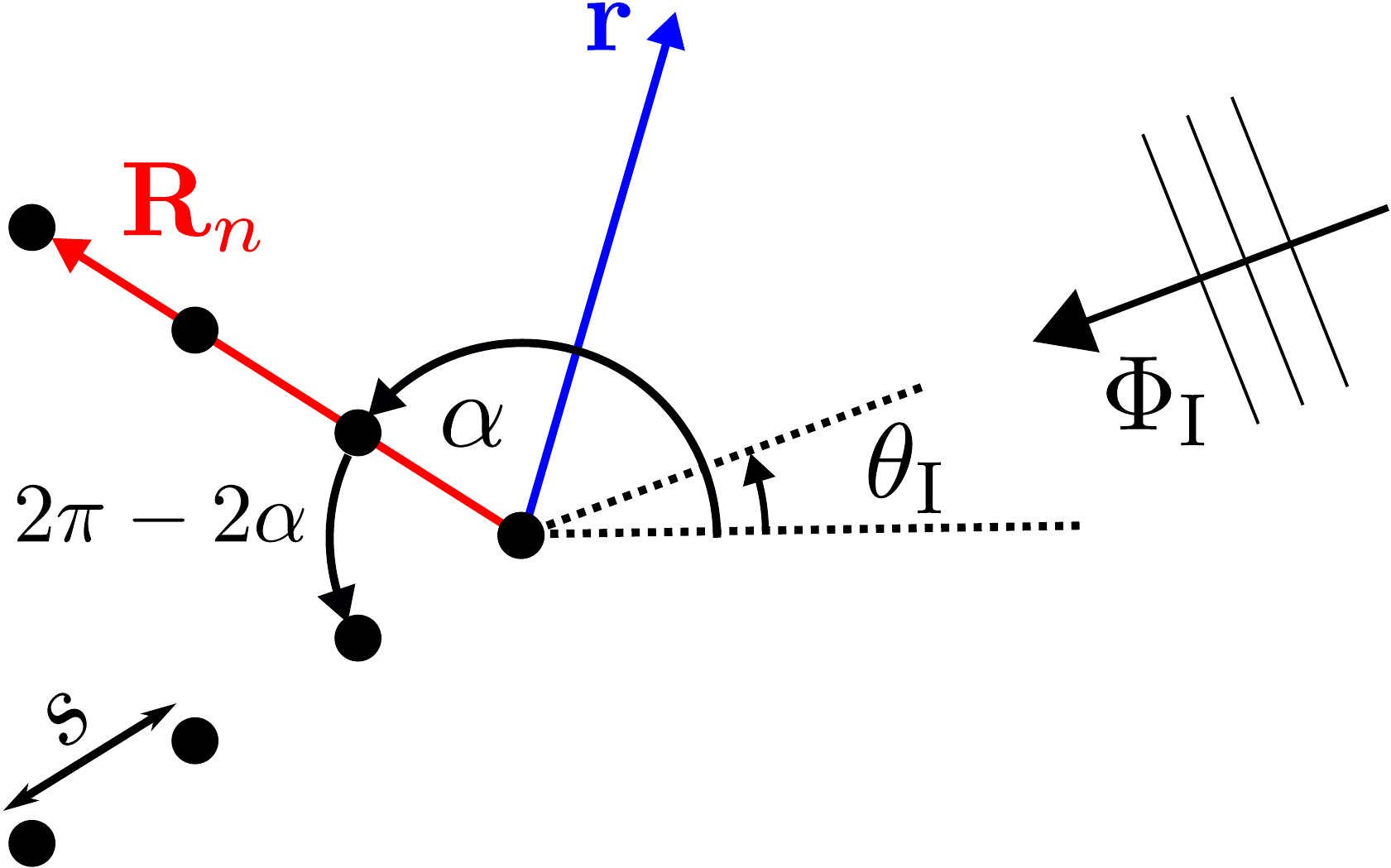}
\caption{Diagram of the point scatterer wedge with scatterers located at $\BF{R}_n$ (see \eqref{HW-R_n}), the position vector $\BF{r}$ and the incident wave $\Phi_{\text{I}}$.}
\label{fig:HW-diagram}
\end{figure}

Now that we have discussed the infinite and semi-infinite array problems, we proceed to the wedge problem (see Figure \ref{fig:HW-diagram}). Here we consider each wedge face to be formed from a semi-infinite array of cylinder scatterers. The positions of these scatterers are given by,
\begin{align}\label{HW-R_n}
\BF{R}_n=&|n|\begin{cases}\BF{s}_+,&n\geq0,\\\BF{s}_-,&n<0,\end{cases}
=|n|s\begin{cases}\cos(\alpha)\BF{\hat{x}}+\sin(\alpha)\BF{\hat{y}},&n\geq0,\\
\cos(\alpha)\BF{\hat{x}}-\sin(\alpha)\BF{\hat{y}},&n<0,\end{cases}
\end{align}
for $n\in\ZZ$. To prevent the cylinders overlapping, this configuration requires $a<s/2$ (as previously stated in Section \ref{IG}) as well as the additional condition $\sin(\alpha)>a/s$. Aside from $\BF{R}_n$, we will use the same definitions and assumptions as Section \ref{IG}. 
For the total field $\Phi$, we pose a solution of the form, 
\begin{align}\label{HW-gensol}
\Phi(\BF{r})=\Phi_{\text{I}}+\sum_{n=0}^\infty \left[A_nH^{(1)}_0(k|\BF{r}-\BF{R}_n|)\right]+\sum_{n=-\infty}^{-1} \left[B_nH^{(1)}_0(k|\BF{r}-\BF{R}_n|)\right],
\end{align}
where $A_n$ and $B_n$ are the coefficients for scatterers on the upper and lower arrays $n\geq0$ and $n<0$ respectively. Because of this restriction, we set $A_n=0$ for all $n<0$ and $B_n=0$ when $n\geq0$. 

\begin{figure}[h!]\centering
\includegraphics[width=0.9\textwidth]{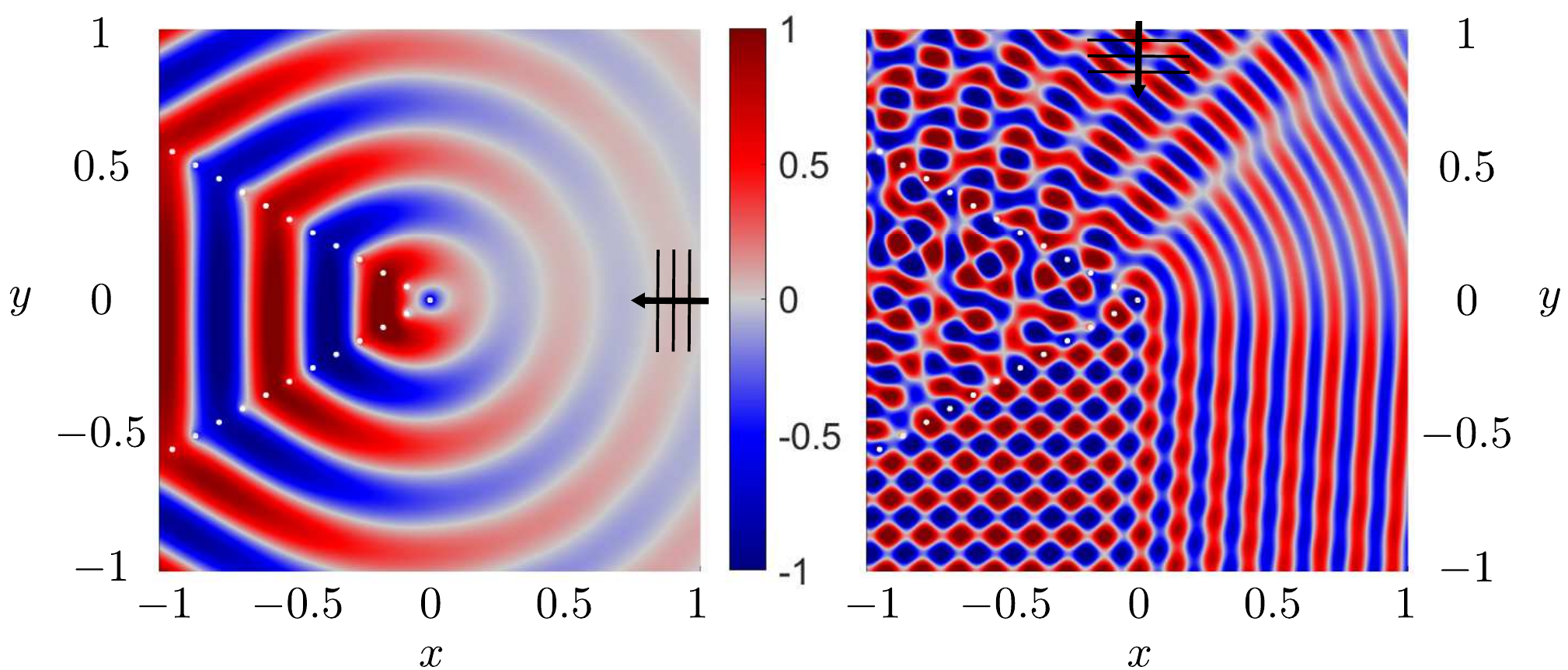}
\caption{The semi-infinite array approximation \eqref{HW-AB-initial} with the wedge defined by $\alpha=\frac{5\pi}{6}$, $s=0.1$ and $a=0.01$. Here we plot the real part of the scattered wave with the left (resp. right) side having an incident wave defined by $k=5\pi$ and $\theta_{\textrm{I}}=0$ (resp. $k=15\pi$ and $\theta_{\textrm{I}}=\frac{\pi}{2}$).}
\label{fig:HW-SIG-approx}
\end{figure}

Initially, we shall approximate the scattering coefficients $A_n$ and $B_n$ by assuming that the two arrays are isolated from each other. This gives us two separate semi-infinite array problems which we can relatively easily find an exact solution for, 
\begin{align}\label{HW-AB-initial}
A_n\approx A^{(0)}_n&= -\frac{e^{-iksn\cos(\theta_{\text{I}}-\alpha)}}{K^+(e^{-iks\cos(\theta_{\text{I}}-\alpha)})}\sum_{m=0}^n\lambda_me^{iksm\cos(\theta_{\text{I}}-\alpha)}\equiv A^{\text{(s-inf)}}_n,\\
\nonumber B_{-n}\approx B^{(0)}_{-n}&= -\frac{e^{-iksn\cos(\theta_{\text{I}}+\alpha)}}{K^+(e^{-iks\cos(\theta_{\text{I}}+\alpha)})}\sum_{m=0}^{n-1}\lambda_me^{iksm\cos(\theta_{\text{I}}+\alpha)}.
\end{align}
The isolation approximation given by \eqref{HW-AB-initial} allows us to get an early glimpse into the behaviour of the scattering coefficients as $n$ becomes large and will be quite useful later in this article as a preliminary approximation. Figure \ref{fig:HW-SIG-approx} illustrates an example of this semi-infinite array approximation including scattering coefficients up to and including $A_{1000}$ and $B_{-1000}$. 

Following the same reasoning as in Section \ref{IG}, starting from \eqref{IG-gensol-phis} and applying the boundary conditions, we can form two coupled infinite systems of equations depending on the value of the index $m$ to solve for the scattering coefficients $A_m$ and $B_m$. If the index is positive ($m\geq0$), then we solve a system of equations for the top array coefficient $A_m$ where the forcing has dependency on the incident wave and the bottom array coefficients, 
\begin{align}\label{HW-Am-system}
A_mH^{(1)}_0(ka)&+\sum_{\substack{n=0\\n\neq m}}^\infty \left[A_nH^{(1)}_0(ks|m-n|)\right]\\
\nonumber &=-\sum_{n=-\infty}^{-1}\left[B_nH^{(1)}_0\left(ks\Lambda_\alpha(m,-n)\right)\right]-e^{-iksm\cos(\theta_{\text{I}}-\alpha)},\quad m\geq0,\\
\label{HW-Lambda-alpha}\textrm{where}\ \ \Lambda_\alpha(m&,n)=(m^2+n^2-2mn\cos(2\alpha))^{\frac{1}{2}}.
\end{align}
If the index is strictly negative ($m<0$), then we solve for $B_m$ in terms of the coefficients $A_n$ and $B_n$,
\begin{align}\label{HW-Bm-system}
B_mH^{(1)}_0(ka)&+\sum_{\substack{n=-\infty\\n\neq m}}^{-1} \left[B_nH^{(1)}_0(ks|m-n|)\right]\\
\nonumber &=-\sum_{n=0}^\infty \left[A_nH^{(1)}_0\left(ks\Lambda_\alpha(-m,n)\right)\right]-e^{iksm\cos(\theta_{\text{I}}+\alpha)},\quad m<0.
\end{align}
From left to right, each of the terms in these two systems have their own physical meaning: the scattered wave from the $m$th cylinder, the interaction with the rest of the cylinders in the associated array, the forcing due to the other array and the forcing due to the incident wave. Note that if we neglect the third terms in these systems (the forcing due to the other array), we are left with the semi-infinite array problems resulting from the isolation assumption (see \eqref{HW-AB-initial}). Solving these two systems will give us $A_m$ in terms of $B_m$ and vice versa.


\subsection{Discrete Wiener--Hopf technique on top array scattering coefficients}\label{HW-A}
We can solve the top array system of equations \eqref{HW-Am-system} using the discrete Wiener--Hopf technique in the same way as the semi-infinite array problem. To start, we need to extend \eqref{HW-Am-system} using some unknown coefficients $F_m$ and recalling $A_m=0$ for $m<0$,
\begin{align}\label{HW-A-begin-system}
A_mH^{(1)}_0(ka)&+\sum_{\substack{n=-\infty\\n\neq m}}^\infty\!\left[A_nH^{(1)}_0\!(ks|m-n|)\right]\\
\nonumber &=\begin{dcases}
-\sum_{n=-\infty}^{-1}\!\left[B_nH^{(1)}_0\!\left(ks\Lambda_\alpha(m,-n)\right)\right]\!-e^{-iksm\cos(\theta_{\text{I}}-\alpha)},&m\geq0,\\
F_m,&m<0.\end{dcases}
\end{align}
Here the $A_m$ scattering coefficients (which are zero for all $m<0$, as stated previously) are the unknowns to find and the $B_m$ coefficients (which are zero for all $m\geq0$) are assumed to be known. To obtain the Wiener--Hopf equation, we apply the same Z-transform as in the semi-infinite array problem, 
\begin{align}\label{HW-A-Z-transform}
G(z)=\sum_{m=-\infty}^\infty G_mz^m,\quad G_m=\frac{1}{2\pi i}\oint_{|z|=1} G(z)z^{-m-1}\text{d}z,
\end{align}
to \eqref{HW-A-begin-system} to obtain the Wiener--Hopf equation
\begin{align}\label{HW-A-WHE}
K(z)A^+(z)=&F_{\text{pole}}^++F_{B}^+(z)+F^-(z),
\end{align}
where $A^+(z)$ is the transform of the unknown scattering coefficients given by
\begin{align}\label{HW-A+def}
A^+(z)=\sum_{m=-\infty}^{\infty} A_mz^m=\sum_{m=0}^{\infty} A_mz^m.
\end{align}
Like the semi-infinite array problem, it is useful to assume that $k$ has a small positive imaginary part to help with the convergence of the Z-transform. We also reuse the definitions of $\Omega^\pm$ (see \eqref{SIG-regions}), to indicate the analytic regions of a function with a $+$ or $-$ superscript.

Of the other components in this Wiener--Hopf equation, the kernel $K(z)$ has the exact same definition and properties as \eqref{SIG-K(z)} including the important identity $K(z)=K(1/z)$ and the singular points $z=e^{\pm iks}$. The three forcing terms on the right-hand side are defined by
\begin{align}
\label{HW-A-Fpole}F_{\text{pole}}^+(z)&=-\sum_{m=0}^\infty e^{-iksm\cos(\theta_{\text{I}}-\alpha)}z^m=\frac{1}{ze^{-iks\cos(\theta_{\text{I}}-\alpha)}-1},\\
\label{HW-A-FB(z)}F_{B}^+(z)&=-\sum_{m=0}^\infty \sum_{n=1}^{\infty}B_{-n}z^mH^{(1)}_0\left(ks\Lambda_\alpha(m,n)\right),\\
\label{HW-A-F-(z)}F^-(z)&=\sum_{m=-\infty}^{-1}F_mz^m.
\end{align}
Here, $F^-(z)=O\left(\frac{1}{z}\right)$ as $|z|\rightarrow\infty$ by design.

We proceed with the Wiener--Hopf technique by factorising $K(z)$ in the exact same way as Section \ref{IG}. This means writing $K(z)=K^+(z)K^-(z)$ in such a way that the two factors satisfy $K^-(1/z)=K^+(z)$ and are defined by Cauchy's integral formula \eqref{SIG-K^+(z)}-\eqref{SIG-K0}. Now let's divide \eqref{HW-A-WHE} by $K^-(z)$,
\begin{align}\label{HW-A-WHE-fac}
K^+(z)A^+(z)=&\frac{F_{\text{pole}}^+(z)}{K^-(z)}+\frac{F_{B}^+(z)}{K^-(z)}+\frac{F^-(z)}{K^-(z)}.
\end{align}
Next, we sum-split each of the following functions,
$$\frac{F_{\text{pole}}^+(z)}{K^-(z)}\ \ \text{and}\ \ \frac{F_{B}^+(z)}{K^-(z)},$$
into two parts that are analytic on $\Omega^\pm$. For the first function, sum-splitting is done by pole removal,
\begin{align}\label{HW-A-pole-sumsplit}
\frac{F_{\text{pole}}^+(z)}{K^-(z)}=&\ \frac{F_{\text{pole}}^+(z)}{K^-(e^{iks\cos(\theta_{\text{I}}-\alpha)})}
+F_{\text{pole}}^+(z)\left[\frac{1}{K^-(z)}-\frac{1}{K^-(e^{iks\cos(\theta_{\text{I}}-\alpha)})}\right].
\end{align}
Here the first term of the right-hand side is analytic on $\Omega^+$ and the second term is analytic on $\Omega^-$. The second function $\frac{F_{B}^+(z)}{K^-(z)}$ is sum-split by rewriting it as a Laurent's series because it is analytic on a thin annulus region of $z$ containing $\Omega^+\cap\Omega^-$.
\begin{align}\label{HW-A-F+Laurent}
\frac{F_{B}^+(z)}{K^-(z)}=D(z)=\sum_{n=-\infty}^\infty D_nz^n,\ \ \text{where}\ \ D_n=\frac{1}{2\pi i}\int_{C^-}D(z)z^{-n-1}\text{d}z.
\end{align}
With this Laurent series, we can easily perform the sum-split, $D(z)=D^+(z)+D^-(z)$,
\begin{align}\label{HW-A-F+sumsplit}
D^+(z)=\sum_{n=0}^\infty D_nz^n,\ \ D^-(z)=\sum_{n=-\infty}^{-1} D_nz^n.
\end{align}
This results in the final Wiener--Hopf equation,
\begin{align}\label{HW-A-WHE-final-LHS}
K^+(z)A^+(z)&-\frac{F_{\text{pole}}^+(z)}{K^-(e^{iks\cos(\theta_{\text{I}}-\alpha)})}-D^+(z)\\
\label{HW-A-WHE-final-RHS}=\ &F_{\text{pole}}^+(z)\left[\frac{1}{K^-(z)}-\frac{1}{K^-(e^{iks\cos(\theta_{\text{I}}-\alpha)})}\right]+D^-(z)+\frac{F^-(z)}{K^-(z)},
\end{align}
where the left-hand (resp. right-hand) side is analytic on $\Omega^+$ (resp. $\Omega^-$). The two sides of this Wiener--Hopf equation are used to construct a function that is entire on the $z$ complex plane,
\begin{align}
\Psi(z)=\begin{cases}
\eqref{HW-A-WHE-final-LHS} & z\in\Omega^+,\\
\eqref{HW-A-WHE-final-RHS} & z\in\Omega^-,\\
\eqref{HW-A-WHE-final-LHS}=\eqref{HW-A-WHE-final-RHS} & z\in\Omega^+\cap\Omega^-.\\
\end{cases}
\end{align}
Each term on the right-hand side is $O\left(\frac{1}{z}\right)$ as $|z|\rightarrow\infty$, which implies that $\Psi$ is bounded and tends to zero at infinity as well as entire. Therefore, Liouville's theorem ensures that both \eqref{HW-A-WHE-final-LHS} and \eqref{HW-A-WHE-final-RHS} are equivalently zero, and then we obtain the solution for $A^+(z)$,
\begin{align}\label{HW-A-WHE-sol}
A^+(z)=\frac{F_{\text{pole}}^+(z)}{K^-(e^{iks\cos(\theta_{\text{I}}-\alpha)})K^+(z)}+\frac{D^+(z)}{K^+(z)}.
\end{align}

To get the scattering coefficients, we must apply the inverse Z-transform to \eqref{HW-A-WHE-sol},
\begin{align}\label{HW-A-integral}
A_m=&\ \frac{e^{iks\cos(\theta_{\text{I}}-\alpha)}}{2\pi iK^-(e^{iks\cos(\theta_{\text{I}}-\alpha)})}
\oint_{C^-} \frac{z^{-m-1}}{K^+(z)(z-e^{iks\cos(\theta_{\text{I}}-\alpha)})}\text{d}z\\
\nonumber &+\frac{1}{2\pi i}\oint_{C^-} \frac{D^+(z)}{K^+(z)}z^{-m-1}\text{d}z.
\end{align}
The integration contour $C^-$ takes the path that runs around the unit circle with $z=0$ being the only potential pole inside. The left integral in \eqref{HW-A-integral} is the solution for the top array in isolation $A^{(0)}_m$, which is equivalent to the semi-infinite array solution $A^{\text{(s-inf)}}_m$ (see \eqref{SIG-A(z)} and \eqref{SIG-Z-inverse}), that is
\begin{align}\label{HW-A-integral-left}
\frac{e^{iks\cos(\theta_{\text{I}}-\alpha)}}{2\pi iK^-(e^{iks\cos(\theta_{\text{I}}-\alpha)})}
\oint_{C^-} \frac{z^{-m-1}}{K^+(z)(z-e^{iks\cos(\theta_{\text{I}}-\alpha)})}\text{d}z
=A^{(0)}_m.
\end{align}
To evaluate the right integral in \eqref{HW-A-integral}, we should recall the Laurent series, \eqref{SIG-1/K-series} and \eqref{HW-A-F+sumsplit},
\begin{align}\label{HW-A-integral-right}
\frac{1}{2\pi i}\oint_{C^-} \frac{D^+(z)}{K^+(z)}z^{-m-1}\text{d}z=\sum_{p=0}^{\infty}\sum_{n=0}^{\infty}\frac{\lambda_pD_n}{2\pi i}\oint_{C^-} z^{p+n-m-1}\text{d}z.
\end{align}
The value of this integral is non-zero only when $p+n-m=0$. This implies that the only non-zero value of the $p$-summation is $p=m-n\geq0$ which implies that,
\begin{align}\label{HW-A-integral-right-sol}
\frac{1}{2\pi i}\oint_{C^-} \frac{D^+(z)}{K^+(z)}z^{-m-1}\text{d}z=\sum_{n=0}^{m}\lambda_{m-n} D_n.
\end{align}
The coefficients $D_n$ are given by \eqref{HW-A-F+Laurent},
\begin{align}\label{HW-A-D-int}
D_n&=-\sum_{q=1}^\infty\sum_{p=0}^\infty\sum_{l=0}^\infty \lambda_p B_{-q} H^{(1)}_0\left(ks\Lambda_\alpha(q,l)\right)
\frac{1}{2\pi i}\int_{C^-}z^{l-n-p-1}\text{d}z,
\end{align}
where the integral is non-zero only when $l-n-p=0$ which implies that,
\begin{align}\label{HW-A-D}
D_n&=-\sum_{q=1}^\infty\sum_{p=0}^\infty \lambda_p B_{-q} H^{(1)}_0\left(ks\Lambda_\alpha(q,p+n)\right).
\end{align}
This means that the scattering coefficients $A_m$ for all $m\geq0$ are written as,
\begin{align}\label{HW-A-sum-sol}
A_m=A^{(0)}_m&-\sum_{q=1}^\infty\sum_{p=0}^\infty\sum_{n=0}^{m} \lambda_{m-n}\lambda_p B_{-q} H^{(1)}_0\left(ks\Lambda_\alpha(q,p+n)\right),\quad m\geq0,
\end{align}
in terms of the coefficients $\lambda_m$ given by \eqref{SIG-lambda-def} and the bottom array coefficients $B_{-m}$. Note that setting $m<0$ in \eqref{HW-A-sum-sol} automatically gives $A_m=0$ as required.

\subsection{Discrete Wiener--Hopf technique on bottom array scattering coefficients}\label{HW-B}
The other system of equations \eqref{HW-Bm-system} can also be solved using the discrete Wiener--Hopf technique for semi-infinite arrays. As before, we extend \eqref{HW-Bm-system} using some unknown coefficients $F_m$ and recalling $B_m=0$ for $m\geq0$, 
\begin{align}\label{HW-B-begin-system}
B_mH^{(1)}_0(ka)&+\sum_{\substack{n=-\infty\\n\neq m}}^{\infty}\! B_nH^{(1)}_0\!(ks|m-n|)\\
\nonumber &=\begin{dcases}F_m,&m\geq0,\\
-\sum_{n=0}^\infty\! A_nH^{(1)}_0\!\left(ks\Lambda_\alpha(-m,n)\right)\!-e^{iksm\cos(\theta_{\text{I}}+\alpha)},&m<0.
\end{dcases}
\end{align}
Contrary to \eqref{HW-A-begin-system}, the $B_m$ scattering coefficients are the unknowns to find and the $A_m$ coefficients are assumed to be known. We also use a slight variation of the Z-transform to \eqref{HW-A-Z-transform} in this problem, 
\begin{align}\label{HW-B-Z-transform}
G(z)=\sum_{m=-\infty}^\infty G_mz^{-m-1},\quad G_m=\frac{1}{2\pi i}\oint_{C^-} G(z)z^{m}\text{d}z.
\end{align}
This new Z-transform is applied to \eqref{HW-B-begin-system} to produce the Wiener--Hopf equation,
\begin{align}\label{HW-B-WHE}
K(z)B^+(z)=&F_{\text{pole}}^+(z)+F_{A}^+(z)+F^-(z),
\end{align}
where $B^+(z)$ is the transform of the bottom array scattering coefficients,
\begin{align}\label{HW-B+def}
B^+(z)=\sum_{m=-\infty}^{\infty} B_mz^{-m-1}=\sum_{m=-\infty}^{-1} B_mz^{-m-1}=\sum_{m=0}^{\infty} B_{-m-1}z^{m}. 
\end{align}
This problem reuses the $\pm$ superscript notation and the definitions of the analytic regions $\Omega^\pm$, the Wiener--Hopf kernel $K(z)$ \eqref{SIG-K(z)} and the factorisations $K^{\pm}(z)$. The three forcing terms on the right-hand side of the Wiener--Hopf equation are defined as,
\begin{align}
\label{HW-B-Fpole}F_{\text{pole}}^+(z)&=-\sum_{m=-\infty}^{-1} e^{iksm\cos(\theta_{\text{I}}+\alpha)}z^{-m-1}=\frac{1}{z-e^{iks\cos(\theta_{\text{I}}+\alpha)}},\\
\label{HW-B-F+(z)}F_{A}^+(z)&=-\sum_{m=1}^{\infty}\sum_{n=0}^\infty A_nz^{m-1}H^{(1)}_0\left(ks\Lambda_\alpha(m,n)\right),\\
\label{HW-B-F-(z)}F^-(z)&=\sum_{m=0}^{\infty} F_mz^{-m-1}.
\end{align}

We proceed with the discrete Wiener--Hopf technique in the same way as in Section \ref{HW-A} and apply the new inverse Z-transform to determine $B_m$. In the end, the bottom array scattering coefficients are given by,
\begin{align}\label{HW-B-sum-sol}
B_{m}=B^{(0)}_{m}&-\sum_{q=0}^\infty\sum_{p=0}^\infty\sum_{n=1}^{-m}\lambda_{-m-n}\lambda_{p}A_q H^{(1)}_0\left(ks\Lambda_\alpha(q,p+n)\right),\quad m<0,
\end{align}
in terms of the coefficients $\lambda_m$ \eqref{SIG-lambda-def} and the top array coefficients $A_m$. Note that setting $m\geq0$ in \eqref{HW-B-sum-sol} automatically gives $B_m=0$ as required. 

\subsection{Point scatterer wedge iterative scheme}
The two formulae for the scattering coefficients given by \eqref{HW-A-sum-sol} and \eqref{HW-B-sum-sol} are still coupled and are difficult to solve exactly. However, we can use them to form an iterative scheme with the isolated solution $A^{(0)}_m$ and $B^{(0)}_{-m}$ as an initial guess 
\begin{align}\label{HW-AB0}
A^{(0)}_m&=-\frac{e^{-iksm\cos(\theta_{\text{I}}-\alpha)}}{K^+(e^{-iks\cos(\theta_{\text{I}}-\alpha)})}\sum_{n=0}^m\lambda_ne^{iksn\cos(\theta_{\text{I}}-\alpha)},\\
\nonumber B^{(0)}_{-m}&=-\frac{e^{-iksm\cos(\theta_{\text{I}}+\alpha)}}{K^+(e^{-iks\cos(\theta_{\text{I}}+\alpha)})}\sum_{n=0}^{m-1}\lambda_ne^{iksn\cos(\theta_{\text{I}}+\alpha)}.
\end{align}
Let us denote $A^{(j)}_m$ and $B^{(j)}_{-m}$ as the $j^{\text{th}}$ approximations to the scattering coefficients $A_m$ and $B_{-m}$ respectively. For the first iteration, we substitute the initial guess $B^{(0)}_{-m}$ into \eqref{HW-A-sum-sol} to get the first iteration of the top array scattering coefficients $A^{(1)}_m$,
\begin{align*}
A^{(1)}_m=A^{(0)}_m-\sum_{q=1}^\infty\sum_{p=0}^\infty\sum_{n=0}^{m}\!\lambda_{m-n}\lambda_p B^{(0)}_{-q} H^{(1)}_0\!\left(ks\Lambda_\alpha(q,p+n)\right).
\end{align*}
Next, we substitute $A^{(1)}_m$ into \eqref{HW-B-sum-sol} to get the first iteration of the bottom array coefficients $B^{(1)}_{-m}$,
\begin{align*}
B^{(1)}_{-m}=B^{(0)}_{-m}-\sum_{q=0}^\infty\sum_{p=0}^\infty\sum_{n=1}^{m}\!\lambda_{m-n}\lambda_{p}A^{(1)}_q H^{(1)}_0\!\left(ks\Lambda_\alpha(q,p+n)\right).
\end{align*}
We can repeat this process as many times as needed to get the $j^{\text{th}}$ approximations, $A^{(j)}_m$ and $B^{(j)}_{-m}$ with the following general equations,
\begin{align*}
A^{(j)}_m&=A^{(0)}_m-\sum_{q=1}^\infty\sum_{p=0}^\infty\sum_{n=0}^{m}\!\lambda_{m-n}\lambda_p B^{(j-1)}_{-q} H^{(1)}_0\!\left(ks\Lambda_\alpha(q,p+n)\right),\\
B^{(j)}_{-m}&=B^{(0)}_{-m}-\sum_{q=0}^\infty\sum_{p=0}^\infty\sum_{n=1}^{m}\!\lambda_{m-n}\lambda_{p}A^{(j)}_q H^{(1)}_0\!\left(ks\Lambda_\alpha(q,p+n)\right).
\end{align*}

Consider the two infinite vectors $\BF{A}^{(j)}$ and $\BF{B}^{(j)}$ which include the $j^{\textrm{th}}$  iteration of the scattering coefficients from $A^{(j)}_0$ to $A^{(j)}_\infty$ and $B^{(j)}_{-1}$ to $B^{(j)}_{-\infty}$ respectively. We can write the iterative scheme as the following infinite matrix system of equations.
\begin{align}\label{HW-AB-matrix-system}
\nonumber \BF{A}^{(j)}&=\BF{A}^{(0)}-\mathcal{M}^{(B)}\BF{B}^{(j-1)},\\
\BF{B}^{(j)}&=\BF{B}^{(0)}-\mathcal{M}^{(A)}\BF{A}^{(j)},
\end{align}
where the matrices $\mathcal{M}^{(A)}$ and $\mathcal{M}^{(B)}$ have the following entries,\\
\begin{align*}
\nonumber \mathcal{M}^{(A)}_{mq}&=\sum_{p=0}^\infty\sum_{n=1}^{m} \lambda_{m-n}\lambda_pH^{(1)}_0\left(ks\Lambda_\alpha(q,p+n)\right),\ \ (m>0,\ q\geq0),\\
\mathcal{M}^{(B)}_{mq}&=\sum_{p=0}^\infty\sum_{n=0}^{m} \lambda_{m-n}\lambda_pH^{(1)}_0\left(ks\Lambda_\alpha(q,p+n)\right),\ \ (m\geq0,\ q>0).
\end{align*}

In practice, we need to truncate the infinite sums at a sufficiently high point. For example, say that we truncate both sums at $M$, then we need to evaluate $M+1$ $\lambda_n$ coefficients (i.e. $\lambda_0,\ \lambda_1,\ ...\ \lambda_M$), for which we use \eqref{Kernel-lambda-approx}. This gives $\mathcal{M}^{(A)}$ as a $M\times(M+1)$ matrix and $\mathcal{M}^{(B)}$ as a $(M+1)\times M$ matrix. Via \eqref{HW-AB0}, that same list of $\lambda_n$ also gives us $M+1$ $A^{(0)}_m$ and $M$ $B^{(0)}_{-m}$ coefficients (i.e. $A^{(0)}_0$, ... $A^{(0)}_M$ and $B^{(0)}_{-1}$, ... $B^{(0)}_{-M}$) as an initial guess. With the iteration matrices and the initial guess defined, we then apply \eqref{HW-AB-matrix-system} $j$ times to get $M+1$ $A^{(j)}_m$ and $M$ $B^{(j)}_{-m}$ coefficients at each stage of the iterative scheme. 

\subsection{Convergence of the iterative scheme}
There are several different parameters that can affect the conditions for convergence of the iterative scheme: $\alpha$, $ks$, $ka$ and the sum truncation numbers. To assess the convergence conditions, we define the $j^{\text{th}}$ error vector of the scattering coefficients by $\BF{E}^{(j)}_A=\BF{A}^{(j)}-\BF{A}$ and $\BF{E}^{(j)}_B=\BF{B}^{(j)}-\BF{B}$ which satisfies the following matrix system, 
\begin{align}\label{HW-error-system}
\nonumber \BF{E}^{(j)}_A&=-\mathcal{M}^{(B)}\BF{E}^{(j-1)}_B,\\
\BF{E}^{(j)}_B&=-\mathcal{M}^{(A)}\BF{E}^{(j)}_A.
\end{align}
We can substitute each equation of \eqref{HW-error-system} into the other to obtain the following,
\begin{align}\label{HW-error-system2}
\nonumber \BF{E}^{(j)}_A&=\mathcal{M}^{(B)}\mathcal{M}^{(A)}\BF{E}^{(j-1)}_A,\\
\BF{E}^{(j)}_B&=\mathcal{M}^{(A)}\mathcal{M}^{(B)}\BF{E}^{(j-1)}_B.
\end{align}
It is well known that the error vectors $\BF{E}^{(j)}_A$ and $\BF{E}^{(j)}_B$ will tend to zero if and only if the spectral radius of either iteration matrices $\mathcal{M}^{(A)}\mathcal{M}^{(B)}$ and $\mathcal{M}^{(B)}\mathcal{M}^{(A)}$ is less than 1, which implies that the iterative scheme will converge. 

\begin{figure}[h!]\centering
\includegraphics[width=0.45\textwidth]{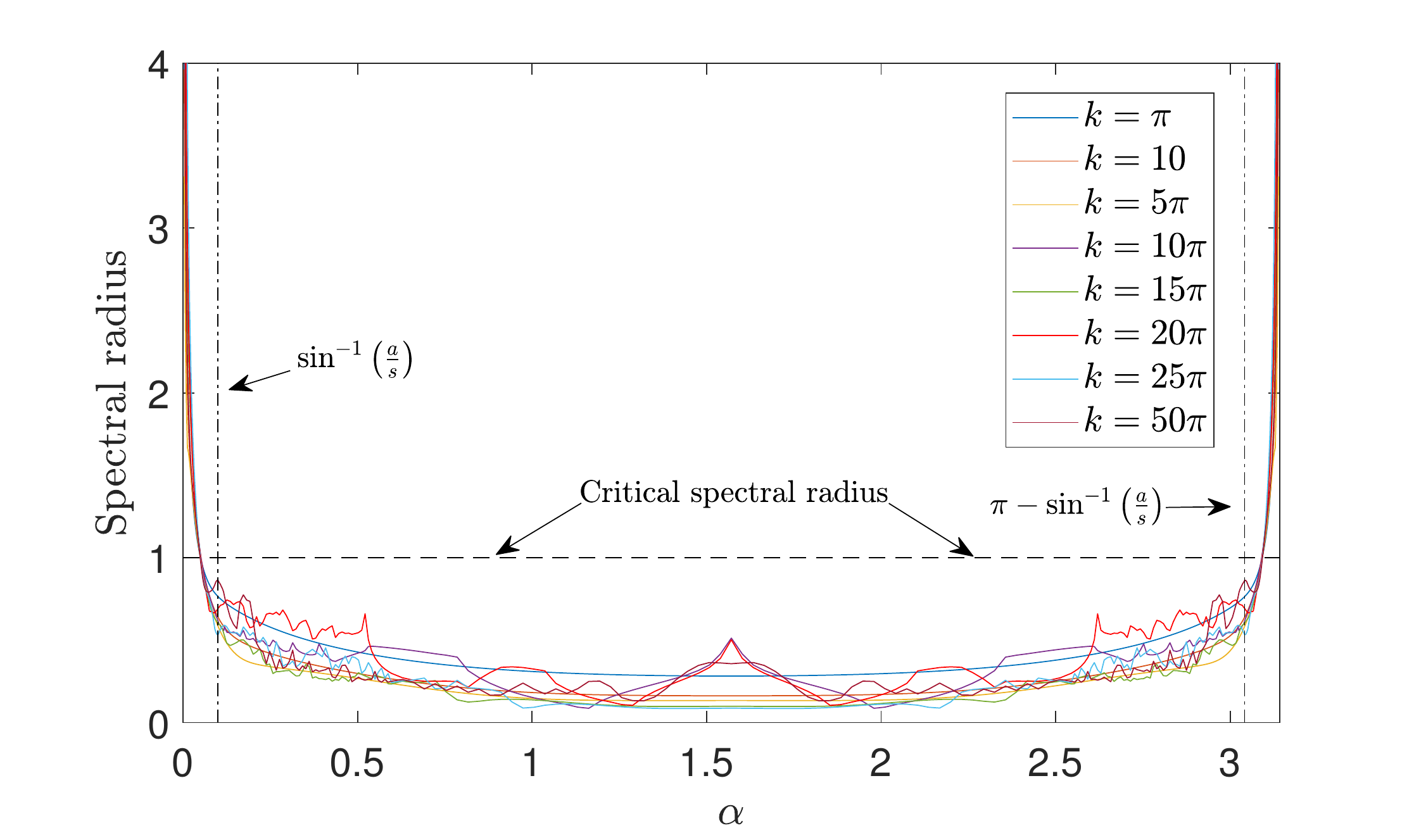}
\includegraphics[width=0.45\textwidth]{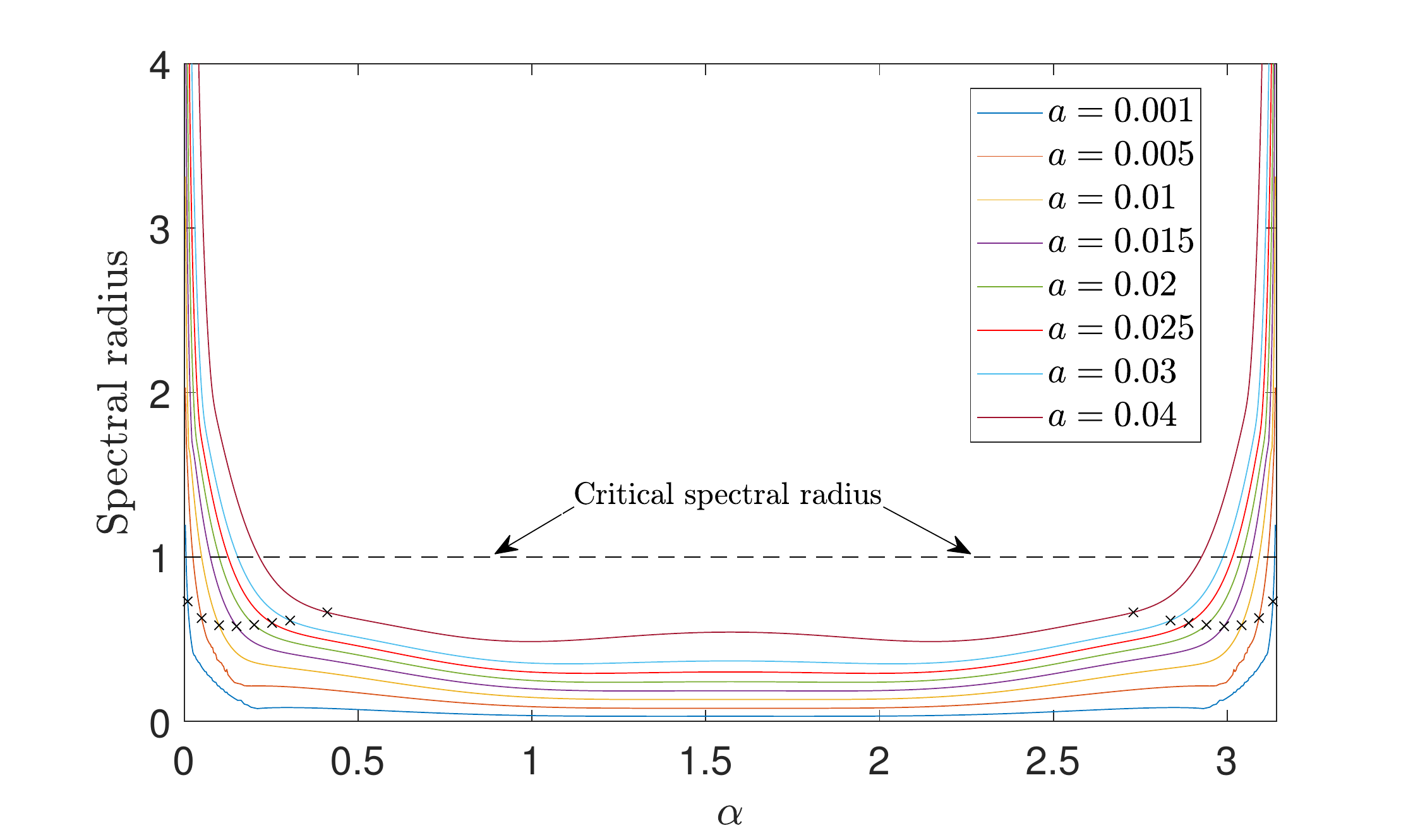}
\caption{Here is the spectral radius $\rho$ of the matrix $\mathcal{M}^{(A)}\mathcal{M}^{(B)}$ with respect to the angle $\alpha$ and summations truncated at $M=1000$. Each plot on the left side has a different value for the wavenumber $k$ while preserving $ka\ll2\pi$. Similarly, the right side varies the cylinder radius $a$ within the constraint $a<s/2$. Note that the horizontal dashed line indicates the critical spectral radius above which the iterative scheme will stop converging. We also mark the edges of the condition $\sin(\alpha)>a/s$ by vertical dashed lines (left) or crosses (right).}
\label{fig:spec_radius}
\end{figure}

We plot the spectral radius of $\mathcal{M}^{(A)}\mathcal{M}^{(B)}$ (and call it $\rho$) in Figure \ref{fig:spec_radius} with respect to the wedge angle $\alpha$. It is important to note that $\rho$ is symmetric about $\alpha=\frac{\pi}{2}$, and although we have plotted $\alpha\in(0,\ \pi)$, the true range of interest here should be $\alpha\in(\sin^{-1}(a/s),\ \pi-\sin^{-1}(a/s))$ which satisfies the condition $\sin(\alpha)>a/s$ (recall this is to prevent overlapping cylinders when $\alpha\approx 0$ or $\pi$). Figure \ref{fig:spec_radius} shows that all the considered cases in this $\alpha$ range will therefore imply that the iterative scheme will converge to a solution.

In producing these graphs, we have observed that $\rho$ converges very quickly with respect to the sum truncation number $M$, so we have deemed it sufficient to set $M=1000$. Figure \ref{fig:spec_radius} (left) illustrates the effect of the wavenumber $k$ on the spectral radius. It shows that the spectral radius is smooth with respect to $\alpha$ when $ks<\pi$, but becomes increasingly jagged as $k$ gets larger (this might be alleviated by increasing $M$). The exception to this is when $ks$ is a multiple of $\pi$, which is a special case that we want to study in a future article. 

Figure \ref{fig:spec_radius} (right) displays how the cylinder radius $a$ affects the convergence. It clearly shows that the iterative scheme will converge faster as $a$ becomes smaller. This trend also continues for larger wavenumbers, however, $a\ll s/2$ will need to be more strictly satisfied to keep $\rho$ below the critical value.

\RED{For higher values of $k$ and $s$, the plots in Figure 4 will appear to be less smooth with respect to $\alpha$. This is due to a change in dominance between the eigenvalues that produce the spectral radius. Apart from some rare exceptions, we can reduce the overall non-smoothness by increasing the value of the sum truncation number $M$, but it should be noted that the computation time is of the order $O(M^3)$.} 

Repeating the same analysis for the other iteration matrix $\mathcal{M}^{(B)}\mathcal{M}^{(A)}$ yields almost identical results (the observed difference is so small that rounding error is the most likely cause). Therefore, we are satisfied that in most cases, the iterative scheme will converge provided that our initial assumptions from Foldy's approximation (namely $a\ll s/2$ and $ka\ll2\pi$) are reasonably satisfied. Additionally, the iterative scheme will perform at its best when $ks<\pi$ as well.

\subsection{Numerical results of the wave field}
Now that we have assessed the spectral radius of the iteration matrices, we have a good idea on what parameters will lead to a successful convergence. However, we have not considered the special case of resonance. Referred to as Wood's anomalies after being observed experimentally in \cite{Wood1902}, this phenomenon occurs when the scattered waves of each unit cell are synchronised and pieced together to form constructive interference as a consequence. For infinite and semi-infinite arrays, the total scattered field usually takes the form of a single plane wave propagating parallel to the array. For the wedge, there are four conditions for resonance (two per array) given by,
\begin{align}\label{HW-res-cond}
\frac{ks}{2\pi}\left(1\mp\cos(\psi-\alpha)\right)\in\ZZ,\ \ \text{and}\ \  \frac{ks}{2\pi}\left(1\mp\cos(\psi+\alpha)\right)\in\ZZ,
\end{align}
where $\psi$ could be the propagation angle of the incident wave or any of the reflected or transmitted wave angles, as determined from geometrical optics. This is a very interesting special case to study and we shall do so in a future article. In the following, we shall choose test cases that are physically viable and non-resonance cases and then, once we are satisfied with the approximations $A^{(j)}_m$ and $B^{(j)}_{-m}$, we substitute them into \eqref{HW-gensol} to get the physical wave fields $\Phi_{\text{S}}$ and $\Phi$. 

\begin{figure}[h!]\centering
\includegraphics[width=0.9\textwidth]{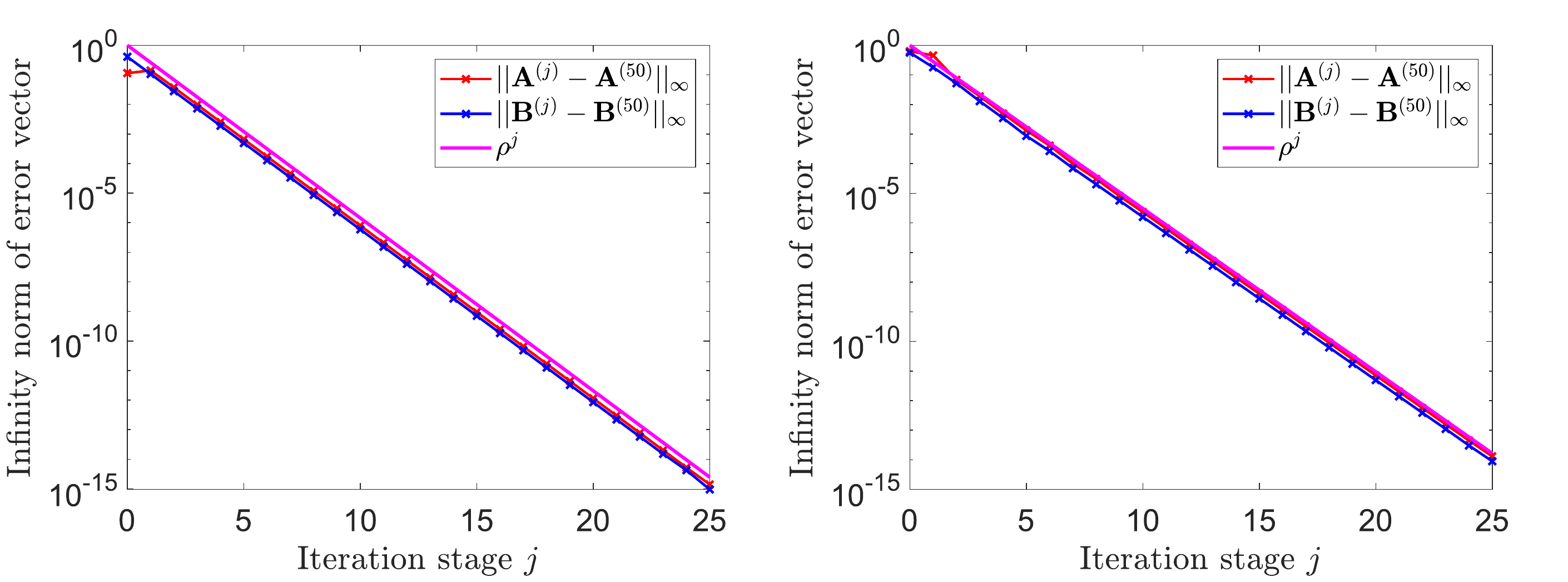}
\caption{We plot the infinity norm error of the scattering coefficients. Here, the {\color{red}red} (resp. {\color{blue}blue}) line displays the error quantity $||\BF{A}^{(j)}-\BF{A}^{(50)}||_{\infty}$ (resp. $||\BF{B}^{(j)}-\BF{B}^{(50)}||_{\infty}$) with respect to iteration stage $j$, where the vectors contain the coefficients up to and including $A^{(j)}_{100}$ (resp. $B^{(j)}_{-100}$) The {\color{magenta}magenta} line is the order of the convergence given by the spectral radius $\rho$. The summations are truncated at $M=1000$, the wedge is defined by $\alpha=\frac{5\pi}{6}$, $s=0.1$ and $a=0.01$, and the left (resp. right) plot has an incident wave with the parameters $k=5\pi$ and $\theta_{\textrm{I}}=0$ (resp. $k=15\pi$ and $\theta_{\textrm{I}}=\frac{\pi}{2}$).}
\label{fig:AB-convergence}
\end{figure}

For two of these test cases, we perform fifty iterations and, in Figure \ref{fig:AB-convergence}, plot the infinity norm of the difference between the $j^{\textrm{th}}$ iteration and the final iteration of the scattering coefficients. In other words, we plot the quantities $||\BF{A}^{(j)}-\BF{A}^{(50)}||_{\infty}$ and $||\BF{B}^{(j)}-\BF{B}^{(50)}||_{\infty}$ in {\color{red}red} and {\color{blue}blue} respectively. Note that the vectors here include coefficients up to and including $A^{(j)}_{100}$ and $B^{(j)}_{-100}$. The {\color{magenta}magenta} line is given by $\rho^j$ where $\rho$ is the spectral radius, which shows that all these errors have the order $O(\rho^j)$. For these test cases, this figure shows that the iteration scheme is becoming accurate to machine precision by the $25^{\textrm{th}}$ iteration.

\begin{figure}[h!]\centering
\includegraphics[width=0.9\textwidth]{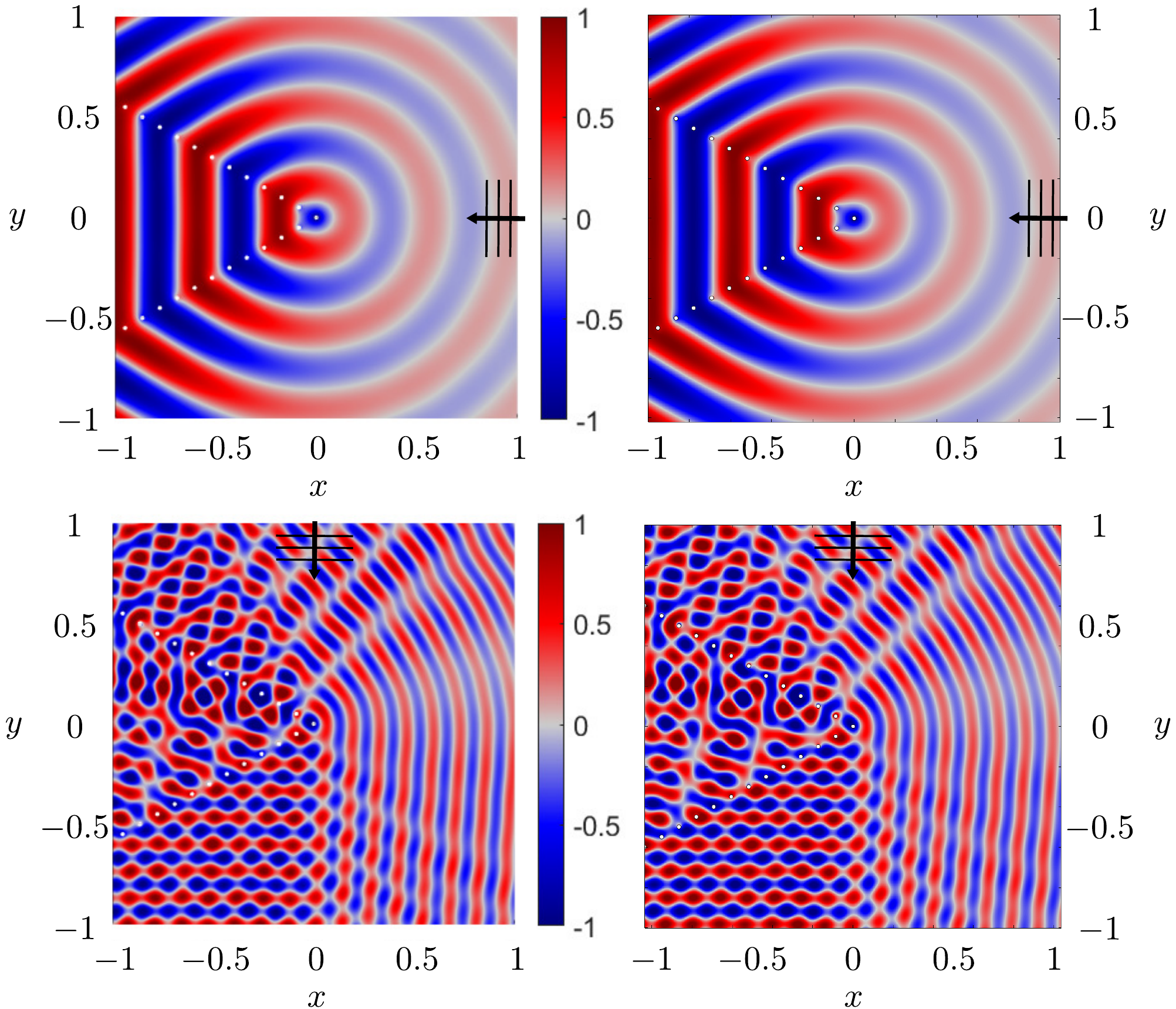}
\caption{Here, we plot the real part of the scattered wave field where the point scatterer wedge is defined by $\alpha=\frac{5\pi}{6}$, $s=0.1$ and $a=0.01$. The left column is using the iterative scheme (twenty-five iterations were performed and summations were truncated at $M=1000$) and the right column plots are COMSOL simulations. The top row has an incident wave with the parameters $k=5\pi$ and $\theta_{\textrm{I}}=0$. The bottom row has the parameters $k=15\pi$ and $\theta_{\textrm{I}}=\frac{\pi}{2}$.}
\label{fig:testcase}
\end{figure}

Figure \ref{fig:testcase} features the scattered wave field's real part of two different test cases using the iterative scheme and compares them to the equivalent finite-element simulations determined through COMSOL. In these figures, we iterated the scattering coefficients twenty-five times and truncated the summations at $M=1000$. This corresponds to considering the interaction of 1000 cylinders on both the top and bottom face as well as the one at the wedge edge. The COMSOL simulations are limited by the number of installed scatterers, as more scatterers requires a larger numerical mesh to find an accurate solution. In these simulations, we created the point scatterer wedge from a total of 61 cylinders and used perfectly-matched layers to prevent incoming radiation. Between the left and right sides of Figure \ref{fig:testcase}, we can clearly see that the agreement is excellent. We have found that in the numerous other test cases we compared, we similarly have excellent agreement.

\begin{figure}[h!]\centering
\includegraphics[width=0.9\textwidth]{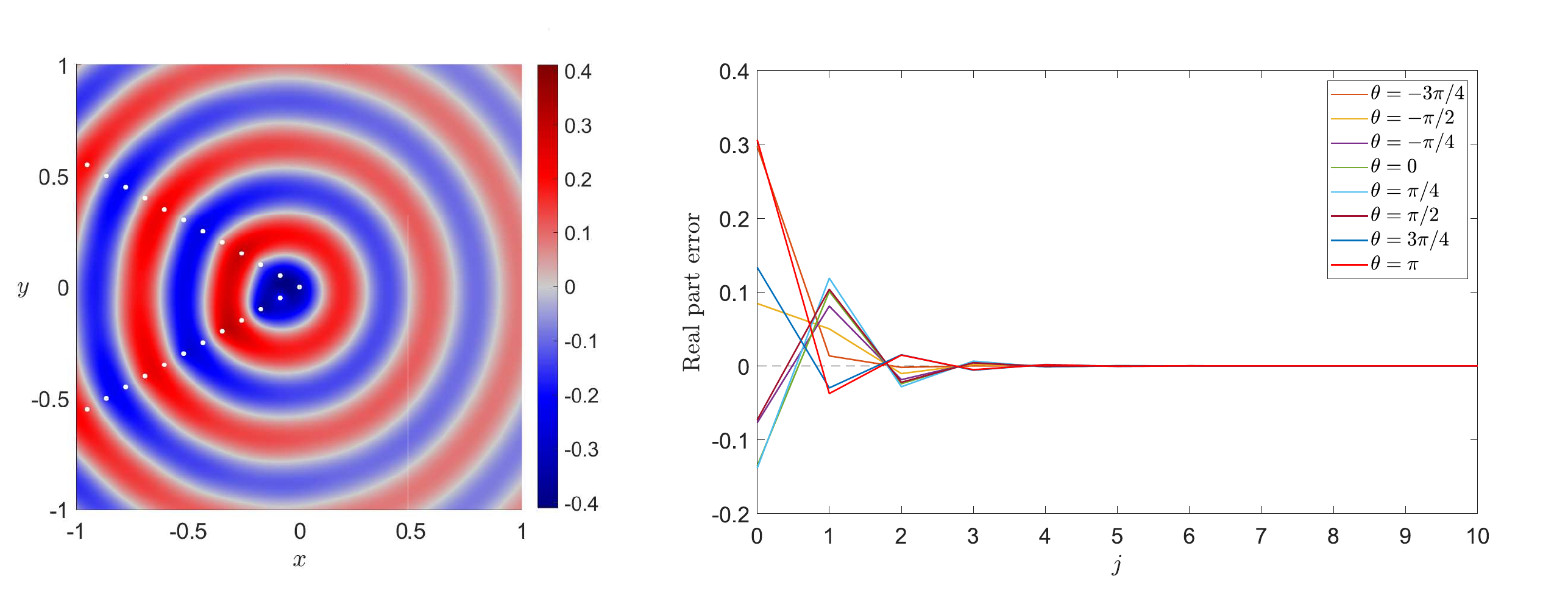}
\caption{The left figure plots the real part difference in the wave field between the $25^{\textrm{th}}$ iteration and the semi-infinite array approximation. The right figure shows the real part error progression of several wave field points. Here, the point scatterer wedge is defined by $\alpha=\frac{5\pi}{6}$, $s=0.1$ and $a=0.01$ and the incident wave has parameters $k=5\pi$ and $\theta_{\textrm{I}}=0$.}
\label{fig:differences}
\end{figure}

To see the evolution of the wave field as the iterative scheme progresses, we plot the difference in the wave field $\Phi$ between the iterative scheme and the semi-infinite array approximation. Specifically, Figure \ref{fig:differences} (left) displays the real part difference between Figure \ref{fig:testcase} (top left) and Figure \ref{fig:HW-SIG-approx} (left) which shows this difference closely resembling a cylindrical wave emanating from the region of the wedge edge. Figure \ref{fig:differences} (right) takes $\RE{\Phi(3s/2,\theta)}$ at several $\theta$ values and illustrates the error between the first ten iterations and the final $25^{\textrm{th}}$ iteration. It shows that, at most points, the convergence is oscillatory. Further plots (not included here) of the absolute error give the convergence order as $O(\rho^j)$, the same behaviour as shown in Figure \ref{fig:AB-convergence}. 


\section{Conclusions and further work}\label{ConFW}
To summarise, we have decomposed the two semi-infinite arrays in the point scatterer wedge into two separate semi-infinite array problems and solved both using the discrete Wiener--Hopf technique. Then, we took this solution and constructed an iterative scheme to approximate the scattering coefficients. Lastly, we assessed the convergence of the iterative scheme via the spectral radius of the iteration matrices and did a side-by-side comparison of some test cases with COMSOL simulations.

We concluded from Figure \ref{fig:spec_radius} that the iterative scheme works very well provided that $\sin(\alpha)>\frac{a}{s}$ is satisfied, we are away from resonance and the assumptions associated with Foldy's approximation ($a\ll\frac{s}{2}$ and $ka\ll 2\pi$) are sufficiently satisfied (in particular, the best convergence occurs when $ks<\pi$). Figure \ref{fig:AB-convergence} showed that the iterative scheme has a convergence order of $O(\rho^j)$ with spectral radius $\rho$ and number of iterations $j$. We also concluded from Figure \ref{fig:testcase} that the iterative scheme agrees closely with the COMSOL comparison and Figure \ref{fig:differences} illustrated how the wave field converges to a solution with respect to the iterations. We see that the first iteration is adding the coupling effect between the two semi-infinite arrays, while each subsequent iteration corrects this coupling effect towards convergence. We also see that the coupling effect is likely to be strongest at the wedge edge and this convergence is most likely exponentially decaying and oscillatory.

Although we mentioned earlier that we deliberately chose test cases that do not result in resonance, this does not mean that we are not interested in the subject. Taking a naive approach, one can say that there are four conditions for resonance which are given by,
\begin{align}\label{ConFW-res-conds}
\frac{ks}{2\pi}(1\mp\cos(\alpha-\theta_{\textrm{I}}))\in\ZZ,\quad\frac{ks}{2\pi}(1\mp\cos(\alpha+\theta_{\textrm{I}}))\in\ZZ,
\end{align}
two for the top array and two for the bottom array, respectively. However, this does not take in account the interactions between the two arrays. For example, could any of the scattered waves (or even resonance waves) from one array cause resonance in the other one? If so, then there will be extra conditions where $\theta_{\textrm{I}}$ in \eqref{ConFW-res-conds} is appropriately replaced for each of these scattered waves. This means that this coupling effect can become convoluted very quickly without adequate ray-tracing techniques to keep track of these scattered waves and their reflections/transmissions. It is very unlikely that the coupling effect influences the nature of the resonance waves themselves, i.e. they will still be in the form of plane waves, but this is unproven.

For the iterative scheme, there are some avenues which can improve the accuracy and speed of numerical computations. These avenues will be especially useful in parameter regimes with slow convergence. One particular avenue is the method of tail-end asymptotics, which could be used to reduce the error of truncations in infinite summations. For example, see Appendix \ref{App-tail} where we applied it to the Wiener-Hopf kernel $K(z)$. This method is also beneficial for calculating both the spectral radius and the scattering coefficients. However, in the latter, it requires prior knowledge of the asymptotics of the scattering coefficients $A_m$ and $B_{-m}$ as $m\rightarrow\infty$ which can be done with the same ray-tracing techniques mentioned earlier for tracking resonance cases.

\RED{Alternatively, this problem can be approached by separating it into distinct symmetric (labelled `even') and anti-symmetric (labelled `odd') sub-problems. In these sub-problems, the incident field is given by,
\begin{align}\label{ConFW-EO-inc}
\Phi_{\text{I}}^{\textrm{(even)}}(r,\theta)=\frac{1}{2}\left(e^{-ikr\cos(\theta-\theta_{\text{I}})}+e^{-ikr\cos(\theta+\theta_{\text{I}})}\right),\\
\Phi_{\text{I}}^{\textrm{(odd)}}(r,\theta)=\frac{1}{2}\left(e^{-ikr\cos(\theta-\theta_{\text{I}})}-e^{-ikr\cos(\theta+\theta_{\text{I}})}\right),
\end{align}
and the scattering coefficients satisfy the relations $B_{-n}^{\textrm{(even)}}=A_n^{\textrm{(even)}}$ and $B_{-n}^{\textrm{(odd)}}=-A_n^{\textrm{(odd)}}$. Each sub-problem leads to a single system of equations to solve:
\begin{align}\label{ConFW-EO-SoE}
\nonumber A^{\textrm{(even)}}_mH^{(1)}_0(ka)&+\sum_{\substack{n=0\\n\neq m}}^\infty \left[A^{\textrm{(even)}}_nH^{(1)}_0(ks|m-n|)\right]\\ \nonumber&+\sum_{n=1}^{\infty}\left[A^{\textrm{(even)}}_nH^{(1)}_0\left(ks\Lambda_\alpha(m,n)\right)\right]+\Phi_{\text{I}}^{\textrm{(even)}}(sm,\alpha)=0,\\
\nonumber A^{\textrm{(odd)}}_mH^{(1)}_0(ka)&+\sum_{\substack{n=0\\n\neq m}}^\infty \left[A^{\textrm{(odd)}}_nH^{(1)}_0(ks|m-n|)\right]\\ &-\sum_{n=1}^{\infty}\left[A^{\textrm{(odd)}}_nH^{(1)}_0\left(ks\Lambda_\alpha(m,n)\right)\right]+\Phi_{\text{I}}^{\textrm{(odd)}}(sm,\alpha)=0,
\end{align}
which are similar to \eqref{HW-Am-system} and \eqref{HW-Bm-system}, but this time the two systems are decoupled. This means that if each system can be solved, then no iteration procedure would be needed. However, each system would only be solvable using the discrete Wiener-Hopf technique provided that we can explicitly write the $Z$-transform of $\sum_{n=1}^{\infty}\left[A_nH^{(1)}_0\left(ks\Lambda_\alpha(m,n)\right)\right]$ in the form of $F(z)A(z)$ where $F(z)$ will become an extra term in the discrete Wiener-Hopf kernel. Since $\Lambda_\alpha(m,n)$ does depend on $\alpha$ and not solely on $m-n$, we cannot use the discrete convolution theorem to do this directly. Hence, at this stage, it is unknown what the explicit formula for $F(z)$ would be. If it was known, then the rest of the procedure should remain exactly the same as the semi-infinite array problem. In the end, the technique would give an exact solution to each sub-problem and then their sum will be the exact solution to the full problem.} 


In this article, we considered the simplest possible unit cell in order to highlight the effects of edges interacting, but certainly many more configurations are possible. Much more complicated models of array scattering with Helmholtz resonators can likewise be considered within the framework of Foldy's approximation \cite{Rodolfo_Ory_asy_20}. We also assumed that the cylinder scatterers were isotropic and hence, restricted the study to sound-soft (or Dirichlet) boundary conditions. However, there are other boundary types that we can consider such as sound-hard (Neumann), impedance or penetrable scatterers for example. In these cases, the isotropic assumption made for Foldy's approximation will not be sufficient. Instead, we will need to at least include double pole terms in the multi-pole expansion which are represented by the Hankel functions $H^{(1)}_{1}(z)$ (see section 8.3.3 in \cite{pmartin2006}). This means that we would have to replace the general solution \eqref{IG-gensol-phis} with,
\begin{align}\label{ConFW-double-poles}
\Phi_{\text{S}}(\BF{r})&=\sum_{n} \left[A_nH^{(1)}_0(kr_n)+(a_n\cos(\theta_n)+b_n\sin(\theta_n))H^{(1)}_{1}(kr_n)\right],
\end{align}
where $r_n=|\BF{r}-\BF{R}_n|$ and $\theta_n=\arg(\BF{r}-\BF{R}_n)$. In this new form, each local scattered wave features the original monopole term as well as an additional dipole term. \RED{However, for different boundary conditions, we anticipate that additional difficulties will arise. This is due to the fact that infinite and semi-infinite arrays with Neumann boundary conditions can support Rayleigh-Bloch surface waves, as discussed in depth in \cite{LintonPorterThompson2007}.}

Apart from Foldy's approximation, another interesting mathematical approach at low frequency is that of asymptotic homogenisation. In particular, for straight or curved arrays of scatterers, the result of such approach is a continuous boundary with effective jump conditions on both the field and its normal derivative \cite{ChapmanHewettTrefethen2015,HewettHewitt2016,TouboulLombardBellis2020,Touboul_etal_2020}. Such boundary is referred to as a meta-interface in the literature. In our case, the meta-interface will be that of a wedge and will hence include a geometric singularity at the corner. It will be very interesting indeed to develop mathematical and/or numerical methods to tackle this `homogenised' point scatterer wedge problem and compare the results to our existing approach.

\RED{Although we designed this iterative scheme with a specific arrangement of scatterers in mind, it is in principle possible to use it with two semi-infinite arrays that do not share a common edge and/or have different spacing parameters ($s_1$ and $s_2$ say). Let us fix one of the arrays to be the traditional semi-infinite array supported on the positive $x$-axis (with $s=s_1$). Then the form of the function $\Lambda_\alpha(m,n)$ will be determined from the relative positioning of the scatterers between the two arrays and will depend on the angle ($\alpha$), the spacing ($s_2$) and the origin of the second array. The only additional necessary condition here for the iterative method to work, is that no two scatterers overlap or touch in any way (as we did for the wedge in this article). Provided that this condition is met, the methodology of this article should be unchanged. Whether or not the iterative scheme converges will be determined by the resulting iterative matrices and their spectral radii.}

\bibliographystyle{abbrvnat}
\bibliography{HWbib}

\begin{thebibliography}{28}
\providecommand{\natexlab}[1]{#1}
\providecommand{\url}[1]{\texttt{#1}}
\expandafter\ifx\csname urlstyle\endcsname\relax
  \providecommand{\doi}[1]{doi: #1}\else
  \providecommand{\doi}{doi: \begingroup \urlstyle{rm}\Url}\fi

\bibitem[Bennetts et~al.(2017)Bennetts, Peter, and Montiel]{Bennetts_Malte_17}
L.~G. Bennetts, M.~A. Peter, and F.~Montiel.
\newblock Localisation of {R}ayleigh-{B}loch waves and damping of resonant
  loads on arrays of vertical cylinders.
\newblock \emph{J. Fluid Mech.}, 813:\penalty0 508--527, 2017.

\bibitem[Brand\~{a}o and Schnitzer(2020)]{Rodolfo_Ory_asy_20}
R.~Brand\~{a}o and O.~Schnitzer.
\newblock Asymptotic modeling of {H}elmholtz resonators including thermoviscous
  effects.
\newblock \emph{Wave Motion}, 97\penalty0 (102583):\penalty0 1--25, 2020.

\bibitem[Bruno and Fernandez-Lado(2017)]{Bruno_Fernandez17}
O.~P. Bruno and A.~G. Fernandez-Lado.
\newblock Rapidly convergent quasi-periodic {G}reen functions for scattering by
  arrays of cylinders---including {W}ood anomalies.
\newblock \emph{Proc. Roy. Soc. A}, 473\penalty0 (2199):\penalty0 20160802,
  2017.

\bibitem[Camacho et~al.(2019)Camacho, Hibbins, Capolino, and Albani]{CHCA2019}
M.~Camacho, A.~P. Hibbins, F.~Capolino, and M.~Albani.
\newblock {Diffraction by a truncated planar array of dipoles: A Wiener–Hopf
  approach}.
\newblock \emph{Wave Motion}, 89:\penalty0 28--42, 2019.

\bibitem[Chapman et~al.(2015)Chapman, Hewett, and
  Trefethen]{ChapmanHewettTrefethen2015}
S.~J. Chapman, D.~P. Hewett, and L.~N. Trefethen.
\newblock Mathematics of the {F}araday cage.
\newblock \emph{SIAM Rev.}, 57\penalty0 (3):\penalty0 398--417, 2015.

\bibitem[Craster and Guenneau(2013)]{Craster-book}
R.~V. Craster and S.~Guenneau.
\newblock \emph{{Acoustic Metamaterials}}.
\newblock Springer, 2013.

\bibitem[Foldy(1945)]{Foldy1945}
L.~L. Foldy.
\newblock {The multiple scattering of waves. I. General theory of isotropic
  scattering by randomly distributed scatterers}.
\newblock \emph{Phys. Rev.}, 67\penalty0 (3-4):\penalty0 107--119, 1945.

\bibitem[Gradshteyn and Ryzhik(2014)]{TablesISP8th}
I.~S. Gradshteyn and I.~M. Ryzhik.
\newblock \emph{{Table of Integrals, Series, and Products}}.
\newblock Academic Press, 8th edition, 2014.

\bibitem[Hewett and Hewitt(2016)]{HewettHewitt2016}
D.~P. Hewett and I.~J. Hewitt.
\newblock Homogenized boundary conditions and resonance effects in {F}araday
  cages.
\newblock \emph{Proc. Roy. Soc. A}, 472\penalty0 (2189), 2016.

\bibitem[Hills and Karp(1965)]{HillsKarp1965}
N.~L. Hills and S.~N. Karp.
\newblock {Semi-Infinite Diffraction Gratings-I}.
\newblock \emph{Comm. Pure Appl. Math}, XVIII:\penalty0 203--233, 1965.

\bibitem[Kamotski\u{\i} and Nazarov(1999)]{Kamotskii_Nazarov99_I}
I.~V. Kamotski\u{\i} and S.~A. Nazarov.
\newblock Wood's anomalies and surface waves in the problem of scattering by a
  periodic boundary. {I}.
\newblock \emph{Mat. Sb.}, 190\penalty0 (1):\penalty0 109--138, 1999.

\bibitem[Korolkov et~al.(2016)Korolkov, Nazarov, and
  Shanin]{Korolkov_Nazarov16}
A.~I. Korolkov, S.~A. Nazarov, and A.~V. Shanin.
\newblock Stabilizing solutions at thresholds of the continuous spectrum and
  anomalous transmission of waves.
\newblock \emph{ZAMM Z. Angew. Math. Mech.}, 96\penalty0 (10):\penalty0
  1245--1260, 2016.

\bibitem[Linton(1998)]{Linton1998}
C.~M. Linton.
\newblock {The {G}reen's function for the two-dimensional {H}elmholtz equation
  in periodic domains}.
\newblock \emph{J. Eng. Math.}, 33\penalty0 (4):\penalty0 377--402, 1998.

\bibitem[Linton and Martin(2004)]{LintonMartin2004}
C.~M. Linton and P.~A. Martin.
\newblock {Semi-Infinite Arrays of Isotropic Point Scatterers. A Unified
  Approach}.
\newblock \emph{SIAM J. Appl. Math}, 64\penalty0 (3):\penalty0 1035--1056,
  2004.

\bibitem[Linton et~al.(2007)Linton, Porter, and
  Thompson]{LintonPorterThompson2007}
C.~M. Linton, R.~Porter, and I.~Thompson.
\newblock {Scattering by a semi-infinite periodic array and the excitation of
  surface waves}.
\newblock \emph{SIAM J. Appl. Math}, 67\penalty0 (5):\penalty0 1233--1258,
  2007.

\bibitem[Liu and Declercq(2015)]{Liu_Declercq15}
J.~Liu and N.~F. Declercq.
\newblock Investigation of the origin of acoustic wood anomaly.
\newblock \emph{J. Acoust. Soc. Am.}, 138\penalty0 (2):\penalty0 1168--1179,
  2015.

\bibitem[Lynott et~al.(2019)Lynott, Andrew, Abrahams, Simon, Parnell, and
  Assier]{Lynott_19}
G.~M. Lynott, V.~Andrew, I.~D. Abrahams, M.~J. Simon, W.~J. Parnell, and R.~C.
  Assier.
\newblock Acoustic scattering from a one-dimensional array; tail-end
  asymptotics for efficient evaluation of the quasi-periodic green's function.
\newblock \emph{Wave Motion}, 89:\penalty0 232--244, 2019.

\bibitem[Martin(2006)]{pmartin2006}
P.~A. Martin.
\newblock \emph{{Multiple Scattering: Interaction of Time-Harmonic Waves with N
  obstacles}}.
\newblock Cambridge University Press, Cambridge, 2006.

\bibitem[Movchan et~al.(2018)Movchan, Movchan, Jones, and
  Colquitt]{MovchanBook}
A.~B. Movchan, N.~V. Movchan, I.~S. Jones, and D.~J. Colquitt.
\newblock \emph{Mathematical Modelling of Waves in Multi-Scale Structured
  Media}.
\newblock {Chapman \& Hall/CRC}, 2018.

\bibitem[Nazarov(2008)]{Nazarov08}
S.~A. Nazarov.
\newblock The {N}eumann problem in angular domains with periodic and parabolic
  perturbations of the boundary.
\newblock \emph{Tr. Mosk. Mat. Obs.}, 69:\penalty0 182--241, 2008.

\bibitem[Nethercote et~al.(2020)Nethercote, Assier, and Abrahams]{WedgeReview}
M.~A. Nethercote, R.~C. Assier, and I.~D. Abrahams.
\newblock {Analytical methods for perfect wedge diffraction: a review}.
\newblock \emph{Wave Motion}, 93\penalty0 (102479), 2020.

\bibitem[Rayleigh(1907)]{Rayleigh1907}
L.~Rayleigh.
\newblock {On the Dynamical Theory of Gratings}.
\newblock \emph{Proc. Roy. Soc. A}, 79:\penalty0 399--416, 1907.

\bibitem[Scarpetta and Sumbatyan(1996)]{Scarpetta_Sumbatyan96}
E.~Scarpetta and M.~A. Sumbatyan.
\newblock Explicit analytical results for one-mode oblique penetration into a
  periodic array of screens.
\newblock \emph{IMA J. Appl. Math.}, 56\penalty0 (2):\penalty0 109--120, 1996.

\bibitem[Shanin and Korolkov(2017)]{Shanin_Korolkov17}
A.~V. Shanin and A.~I. Korolkov.
\newblock Diffraction of a mode close to its cut-off by a transversal screen in
  a planar waveguide.
\newblock \emph{Wave Motion}, 68:\penalty0 218--241, 2017.

\bibitem[Touboul et~al.(2020{\natexlab{a}})Touboul, Lombard, and
  Bellis]{TouboulLombardBellis2020}
M.~Touboul, B.~Lombard, and C.~Bellis.
\newblock Time-domain simulation of wave propagation across resonant
  meta-interfaces.
\newblock \emph{J. Comput. Phys}, 414\penalty0 (109474), 2020{\natexlab{a}}.

\bibitem[Touboul et~al.(2020{\natexlab{b}})Touboul, Pham, Maurel, Marigo,
  Lombard, and Bellis]{Touboul_etal_2020}
M.~Touboul, K.~Pham, A.~Maurel, J.~Marigo, B.~Lombard, and C.~Bellis.
\newblock Effective resonant model and simulations in the time-domain of wave
  scattering from a periodic row of highly-contrasted inclusions.
\newblock \emph{J. Elast}, 142:\penalty0 53--82, 2020{\natexlab{b}}.

\bibitem[Wegert(2012)]{Wegert2012}
E.~Wegert.
\newblock \emph{{Visual Complex Functions}}.
\newblock Birkhauser Basel, 2012.

\bibitem[Wood(1902)]{Wood1902}
R.~W. Wood.
\newblock {On a remarkable case of uneven distribution of light in a
  diffraction grating spectrum}.
\newblock \emph{Proc. Phys. Soc.}, 18\penalty0 (1):\penalty0 269--275, 1902.

\end{thebibliography}

\begin{appendix}
\section{Asymptotics and numerical methods for the discrete Wiener--Hopf kernel and factors}\label{App-Kernel}
Infinite sums are required to be truncated repeatedly in this article. In this appendix section, we shall discuss two different computational methods to evaluate the discrete Wiener--Hopf kernel. The first is the method of tail-end asymptotics to improve the convergence of the original formula. The second is an alternative fast-convergent formula. Recall the original formula for the discrete Wiener--Hopf kernel \eqref{SIG-K(z)},
\begin{align}\label{Kernel-original}
K(z)=H^{(1)}_0(ka)+\sum_{\ell=1}^\infty (z^\ell+z^{-\ell})H^{(1)}_0(ks\ell).
\end{align}
In the limit as $\IM{k}\rightarrow0$, this formula only converges when $z$ is on the unit circle (not including the branch points $z=e^{\pm iks}$) and even then, the convergence is conditional and very slow. This is because the sum term has the following order as $\ell\rightarrow\infty$,
\begin{align}\label{Kernel-K-term-asymp}
(z^\ell+z^{-\ell})H^{(1)}_0(ks\ell)=O\left((z^\ell+z^{-\ell})\frac{e^{iks\ell}}{\sqrt{\ell}}\right).
\end{align}
Computing the kernel will naturally require us to truncate the sum at $L$ say. This introduces an error term which we call the tail-end. 

\subsection{Method of tail-end asymptotics}\label{App-tail}
Let us define the tail-end sum by,
\begin{align}\label{Kernel-tail-end}
\text{Tail}(z,L)=\sum_{\ell=L}^\infty (z^\ell+z^{-\ell})H^{(1)}_0(ks\ell).
\end{align}
If $L$ was sufficiently large, one could neglect the tail-end altogether, but the slow convergence of \eqref{Kernel-original} would make this impractical. Instead, the tail-end sum is asymptotically approximated as $L\rightarrow\infty$ (provided that $z\neq e^{\pm iks}$) using the same methodology as in \cite{Lynott_19} for the quasi-periodic Green's function. With two corrections, the result of this method of tail-end asymptotics is given by,
\begin{align}\label{Kernel-Tasymp}
\text{Tail}(z,L)&\sim\frac{1-i}{\sqrt{\pi ksL}}\left[T(ze^{iks},L)+T(e^{iks}/z,L)\right],\\
\nonumber\text{where}\ T(z,L)&=\frac{z^L}{1-z}\bigg[1-\frac{1}{2L}\left(\frac{z}{1-z}+\frac{i}{4ks}\right)\bigg]. 
\end{align}
Even just one tail-end correction will drastically improve the convergence of the Hankel summation formula. Specifically $N$ tail-end corrections will imply that, 
\begin{align}\label{Kernel-tail-end-form}
K(z)=&\ H^{(1)}_0(ka)+\sum_{\ell=1}^{L-1} (z^\ell+z^{-\ell})H^{(1)}_0(ks\ell)+\text{Tail}(z,L)\\ \nonumber &+O\left(\left(\frac{z^L}{1-ze^{iks}}+\frac{z^{-L}}{1-e^{iks}/z}\right)\frac{e^{iksL}}{L^{\frac{1}{2}+N}}\right).
\end{align}
To reinforce our statements, Figure \ref{fig:K-comparison} compares the various rates of convergence of \eqref{Kernel-tail-end-form} at $z=1$ with a varying truncation point $L$ and different numbers of tail-end corrections $N$. These plots clearly show that no tail-end corrections will cause very slow convergence while even just one tail-end correction is a significant improvement. Evidently, the method of tail-end asymptotics is a powerful technique for slow-convergent infinite sums like \eqref{Kernel-original}. It is also useful for improving the convergence of other infinite sums in this article including the physical field \eqref{IG-gensol-phis} and the scattering coefficients \eqref{IG-A0-SchSeries}, \eqref{SIG-Am-final}, \eqref{HW-A-sum-sol} and \eqref{HW-B-sum-sol}. However, these cases will require us to determine the asymptotic behaviour of the $\lambda_n$ coefficients as $n\rightarrow\infty$.

\subsection{Fast-convergent formula}
Besides the original formula \eqref{Kernel-original}, we can reuse equations (8.522) and (8.524) from \cite{TablesISP8th} to write the discrete Wiener--Hopf kernel in a form that is faster to converge. Before we do this, it is convenient to map $z$ onto a new complex plane using the following transformation $z=e^{it}$. Under this mapping the $z$ unit circle is transformed onto the real line, inside (resp. outside) the unit circle is mapped onto the upper (resp. lower) half plane and \eqref{Kernel-original} is rewritten as 
\begin{align}\label{Kernel-K(e^it)}
K(e^{it})=H^{(1)}_0(ka)+2\sum_{\ell=1}^\infty \cos(\ell t)H^{(1)}_0(ks\ell).
\end{align}
Applying equations (8.522) and (8.524) to \eqref{Kernel-K(e^it)} gives us the fast-convergent formula,
\begin{align}\label{Kernel-K(z)-fast}
K(e^{it})=&\ H^{(1)}_0(ka)-1-\frac{2i}{\pi}\left(\gamma+\ln\left(\frac{ks}{4\pi}\right)\right)+\frac{2}{\sqrt{(ks)^2-t^2}}\\
\nonumber &+\sum_{\substack{\ell=-\infty\\\ell\neq 0}}^\infty\left(\frac{2}{\sqrt{(ks)^2-(2\pi \ell-t)^2}}+\frac{i}{\pi |\ell|}\right).
\end{align}
where $\gamma=0.5772...$ is the Euler–Mascheroni constant. In this formula, there are two sets of branch cuts starting at the points $t=\pm ks+2\pi\ell$. The ones starting at $t=ks+2\pi\ell$ (resp. $t=-ks+2\pi\ell$) are orientated to positive (resp. negative) imaginary infinity in the $t$ complex plane. All of the square roots in \eqref{Kernel-K(z)-fast} are defined to ensure that $\IM{\sqrt{(ks)^2-(2\pi \ell-t)^2}}\geq0$ when $t$ is real and for all $\ell$. In other words, the branches are selected such that $\left.\sqrt{(ks)^2-(2\pi \ell-t)^2}\right|_{t=2\pi \ell}=ks$ for all $\ell$. 

We can further accelerate the convergence of \eqref{Kernel-K(z)-fast} by using the following asymptotic expansion as $\ell\rightarrow\infty$,
\begin{align}\label{Kernel-sqrt-asymp}
\frac{2}{\sqrt{(ks)^2-(2\pi\ell-t)^2}}+\frac{2}{\sqrt{(ks)^2-(2\pi\ell+t)^2}}+\frac{2i}{\pi\ell}
=\frac{2t^2+(ks)^2}{4i(\pi\ell)^3}+O(\ell^{-5}).
\end{align}
Adding and subtracting the right hand side terms of \eqref{Kernel-sqrt-asymp} into \eqref{Kernel-K(z)-fast} will give us the accelerated formula,
\begin{align}
\nonumber K(e^{it})=&\ H^{(1)}_0(ka)-1-\frac{2i}{\pi}\left(\gamma+\ln\left(\frac{ks}{4\pi}\right)\right)+\frac{2}{\sqrt{(ks)^2-t^2}}+\frac{2t^2+(ks)^2}{4i\pi^3}\zeta(3)\\
\label{Kernel-K(z)-faster} &+\sum_{\ell=1}^\infty\Bigg(\frac{2}{\sqrt{(ks)^2-(2\pi\ell-t)^2}}+\frac{2}{\sqrt{(ks)^2-(2\pi\ell+t)^2}}+\frac{2i}{\pi\ell}
-\frac{2t^2+(ks)^2}{4i(\pi\ell)^3}\Bigg),
\end{align}
where $\zeta(z)$ is the Riemann zeta function. More terms in the expansion \eqref{Kernel-sqrt-asymp} will allow us the accelerate the convergence even further, 

\noindent\resizebox{1\linewidth}{!}{
\begin{minipage}{\linewidth}
\begin{align}
\nonumber K(e^{it})=&\ H^{(1)}_0(ka)-1-\frac{2i}{\pi}\left(\gamma+\ln\left(\frac{ks}{4\pi}\right)\right)+\frac{2}{\sqrt{(ks)^2-t^2}}+\frac{2t^2+(ks)^2}{4i\pi^3}\zeta(3)+\frac{8t^4+24(kst)^2+3(ks)^4}{64i\pi^5}\zeta(5)\\
\label{Kernel-K(z)-fastest} &+\sum_{\ell=1}^\infty\Bigg(\frac{2}{\sqrt{(ks)^2-(2\pi\ell-t)^2}}+\frac{2}{\sqrt{(ks)^2-(2\pi\ell+t)^2}}+\frac{2i}{\pi\ell}
-\frac{2t^2+(ks)^2}{4i(\pi\ell)^3}-\frac{8t^4+24(kst)^2+3(ks)^4}{64i(\pi\ell)^5}\Bigg).
\end{align}
\end{minipage}}

\begin{figure}[h!]\centering
\includegraphics[width=0.9\textwidth]{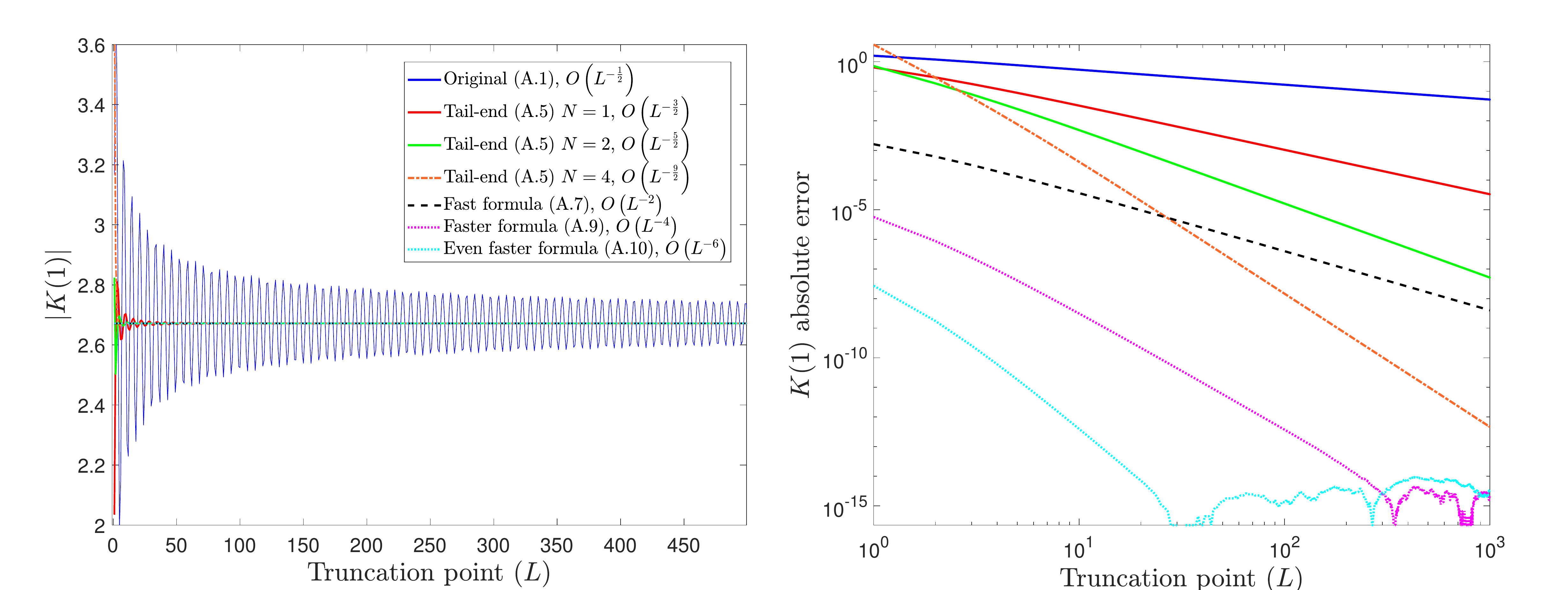}
\caption{These graphs display the progression of the different Wiener--Hopf kernel methods (left) at $z=1$ and their absolute error with a high order reference value (right) with respect to the truncation point $L$. Here we have $ks=1$ and $ka=0.01$. For both figures $K(1)$ is computed using either the original formula \eqref{Kernel-original} (see {\color[rgb]{0,0,1}blue} line), with $N$ additional tail-end corrections \eqref{Kernel-tail-end-form} ($N=1,2$ and $4$ for the {\color[rgb]{1,0,0}red}, {\color[rgb]{0,1,0}green} and {\color[rgb]{1,0.41,0.16}orange} lines respectively) or using the fast-convergent formula (\eqref{Kernel-K(z)-fast} is the black line, \eqref{Kernel-K(z)-faster} is {\color[rgb]{1,0,1}magenta} and \eqref{Kernel-K(z)-fastest} is {\color[rgb]{0,1,1}cyan}). The legend also gives the rate of convergence with respect to $L$ for each line.}
\label{fig:K-comparison}
\end{figure}

If we take \eqref{Kernel-K(z)-fast} with $N$ extra corrections and truncated at $L$, the neglected tail-end of that sum has the asymptotic behaviour $O\left(L^{-(2N+2)}\right)$ as $L\rightarrow\infty$, which is also the rate of convergence. Figure \ref{fig:K-comparison} additionally compares \eqref{Kernel-K(z)-fast}, \eqref{Kernel-K(z)-faster} and \eqref{Kernel-K(z)-fastest} against the method of tail-end asymptotics. It shows that these formula are generally much faster to converge than \eqref{Kernel-tail-end-form}, until it gets to machine precision and then the working precision needs to be increased.  

Another benefit of \eqref{Kernel-K(z)-fast} is the ability to converge for complex values of $t$ (i.e. an annulus about the $z$ unit circle). This ability coupled with the fast convergence, is why we use \eqref{Kernel-K(z)-fastest} for analytic evaluations of the Wiener--Hopf kernel.

\subsection{Rational approximation}
To produce the final results of this article, we will need to create a rational approximation of the discrete Wiener--Hopf kernel. After construction, the rational approximation is computationally very fast, remarkably accurate and is also trivial to factorise analytically for the Wiener--Hopf technique. This rational approximation is given by
\begin{align}\label{Kernel-K_approx-basic}
K(z)\approx \widetilde{K}(z)=\widetilde{K}_1\prod_{\ell=1}^{L}\frac{(z-z_\ell)}{(z-p_\ell)},
\end{align}
where the constant $\widetilde{K}_1$ and the zero and pole locations ($z_\ell$ and $p_\ell$ respectively) are determined for the most accuracy. In MATLAB, the function we use to create this rational approximation is ``AAA" in Chebfun. Here, we input the values of $K(z)$ on a set of points located on the $z$ unit circle using \eqref{Kernel-K(z)-fastest} and then the ``AAA'' function outputs the rational approximation of the form \eqref{Kernel-K_approx-basic} as well as individual lists of zeros ($z_\ell$), poles ($p_\ell$) and associated residues. 

\begin{figure}[h!]\centering
\includegraphics[width=0.9\textwidth]{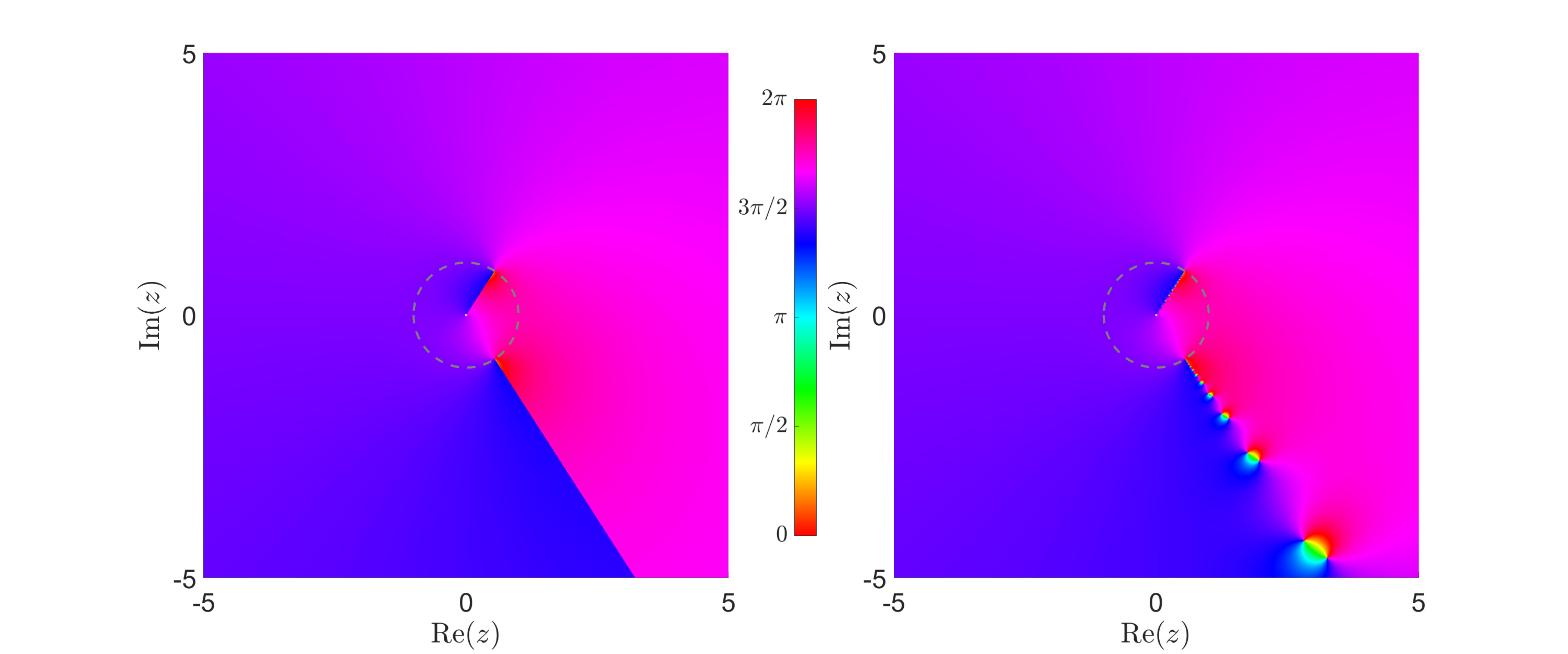}
\caption{Phase portrait comparison between the exact formula for the Wiener--Hopf kernel \eqref{Kernel-K(z)-fastest} (left) and the ``AAA'' approximation \eqref{Kernel-K_approx-basic} (right). Phase portraits indicate the complex argument of a function by using a HSV colour scale. In this example, we have $ka=0.01$ and $ks=1$.}
\label{fig:K-approx-comparison}
\end{figure}

Figure \ref{fig:K-approx-comparison} compares the exact formula \eqref{Kernel-K(z)-fastest} and the ``AAA'' approximation \eqref{Kernel-K_approx-basic} on the $z$ complex plane by plotting the complex argument \cite{Wegert2012}. This comparison shows that the zero and poles of the rational approximation are positioned to best imitate the branch cuts. 

For convenience, the zeros and poles of the approximation are split into two sets depending on their locations being inside or outside the unit circle. Hence, we split and relabel the lists $z_\ell$ and $p_\ell$ with superscripts $-$ and $+$ to indicate their position inside and outside $|z|=1$,
\begin{align}\label{Kernel-K_approx-split}
K(z)\approx \widetilde{K}(z)=\widetilde{K}_1\prod_{\ell=1}^{L^-}\frac{(z-z^-_\ell)}{(z-p^-_\ell)}\prod_{\ell=1}^{L^+}\frac{(z-z^+_\ell)}{(z-p^+_\ell)},
\end{align}
where $L^-+L^+=L$. In these lists, we expect $\arg(z^+_\ell)\approx\arg(p^+_\ell)\approx\arg(e^{-iks})$ and $\arg(z^-_\ell)\approx\arg(p^-_\ell)\approx\arg(e^{iks})$. We also expect that $|z^+_\ell|,|p^+_\ell|>1$ and $|z^-_\ell|,|p^-_\ell|<1$ and then, without loss of generality, the zeros and poles are organised by absolute value from the unit circle, i.e.
\begin{align}\label{Kernel-K_approx-organise}
|z^+_\ell|<|z^+_{\ell+1}|,\ \ |p^+_\ell|<|p^+_{\ell+1}|,\ \ |z^-_\ell|>|z^-_{\ell+1}|,\ \ |p^-_\ell|>|p^-_{\ell+1}|.
\end{align}

\begin{figure}[h!]\centering
\includegraphics[width=0.9\textwidth]{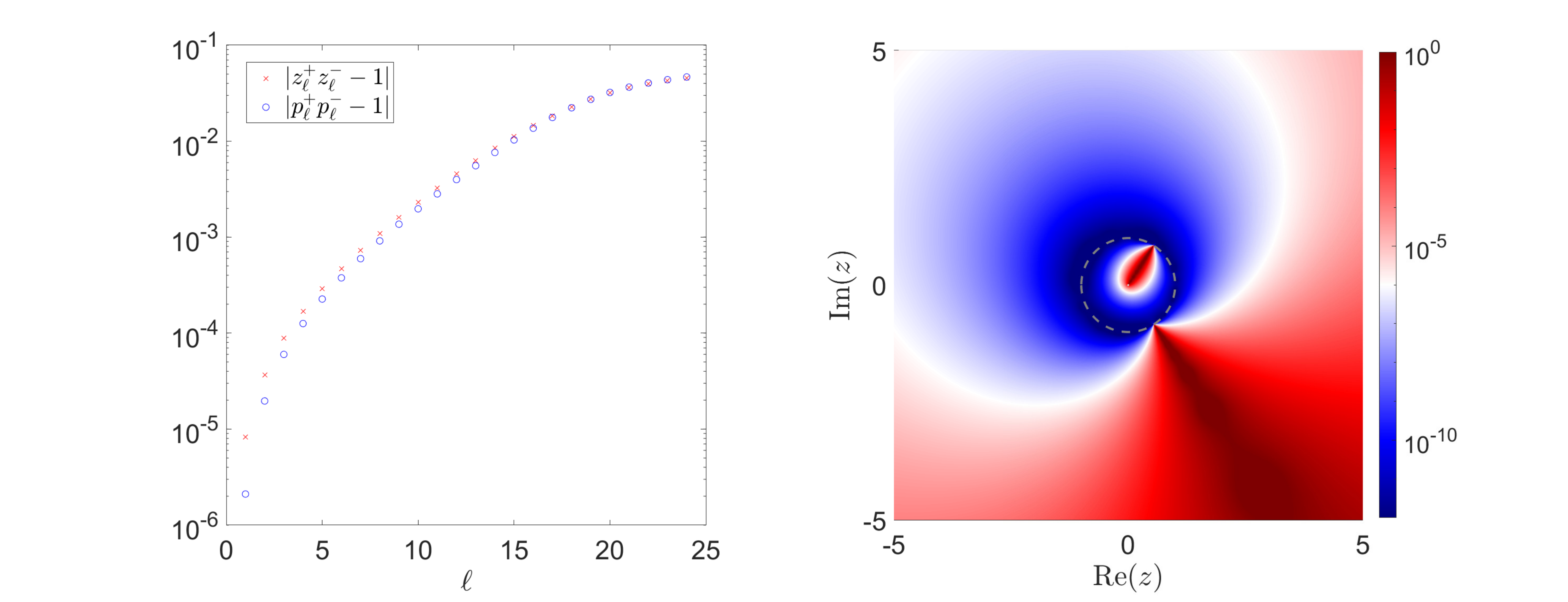}
\caption{On the left side, we plot the absolute value of the two quantities $z^+_\ell z^-_\ell\approx1$ and $p^+_\ell p^-_\ell\approx1$ for all $\ell$. On the right side, we plot the absolute error of the rational approximation $|K(z)-\widetilde{K}(z)|$ where the colour map is scaled logarithmically. In this example, we have $ka=0.01$ and $ks=1$.}
\label{fig:K-reciprocability}
\end{figure}

For this approximation to be effective, it is very important to try to conserve the kernel identity $K(z)=K\left(\frac{1}{z}\right)$ as much as possible. For this, we need the same number of zeros and poles inside the unit circle as outside (meaning $L^-=L^+=L/2$). Ideally, the zeros and poles inside and outside should also be in reciprocal pairs (or at least approximately so), in other words $z^+_\ell z^-_\ell\approx1$ and $p^+_\ell p^-_\ell\approx1$ for all $\ell$ (see Figure \ref{fig:K-reciprocability} (left) for a test case). We are not aware of a way to enforce these properties in the ``AAA" approximation, but we can increase the likelihood by using reciprocal pairs in the $z$ input data and employing checks that restart the approximation with fresh data when the properties are not satisfied enough. Figure \ref{fig:K-reciprocability} (right) plots the absolute error of the approximation. This figure illustrates the high accuracy rating of the rational approximation not only on the unit circle from which it is created, but in its neighbourhood too (the exception being in the vicinity of the branch points and cuts).

The downside to this approximation is that it cannot capture the correct asymptotic behaviour at the branch points, no matter how good the approximation is elsewhere. However, this is to be expected considering the format of the rational approximation. 

When factorising the Wiener--Hopf kernel, we recall how the $+$ (or $-$) factor must be analytic and zero-free inside (or outside) the unit circle. For the approximate kernel, it is clear that the $+$ factor will have zeros at $z^+_\ell$ and poles at $p^+_\ell$. Similarly the $-$ factor will have zeros and poles at $z^-_\ell$ and $p^-_\ell$ respectively. The result of the factorisation is given by,
\begin{align}\label{Kernel-K_approx-fac}
K^+(z)\approx\widetilde{K}^+(z)=\widetilde{K}^+_1\prod_{\ell=1}^{L/2}\frac{z-z^+_\ell}{z-p^+_\ell},\quad 
K^-(z)\approx\widetilde{K}^-(z)=\widetilde{K}^-_1\prod_{\ell=1}^{L/2}\frac{z-z^-_\ell}{z-p^-_\ell},
\end{align}
where the two constants $\widetilde{K}^+_1$ and $\widetilde{K}^-_1$ satisfy $\widetilde{K}^+_1\widetilde{K}^-_1=\widetilde{K}_1$ by design. We chose $\widetilde{K}^+_1$ such that $\widetilde{K}^+(z)=\widetilde{K}^-(1/z)$ is satisfied on at least one point on the unit circle. For example, say $z=1$ (or $z=-1$ to avoid errors due to branch point singularities if $ks\approx 2\pi M$ ($M\in\ZZ$)), which implies that, 
\begin{align}\label{Kernel-K1plus}
\widetilde{K}^+_1=\left(\widetilde{K}_1\prod_{\ell=1}^{L/2}\frac{(1\mp z^-_\ell)(1\mp p^+_\ell)}{(1\mp p^-_\ell)(1\mp z^+_\ell)}\right)^{\frac{1}{2}}
=\left(\widetilde{K}(\pm 1)\right)^{\frac{1}{2}}\prod_{\ell=1}^{L/2}\frac{1\mp p^+_\ell}{1\mp z^+_\ell},
\end{align}
then $\widetilde{K}^-_1=\widetilde{K}_1/\widetilde{K}^+_1$. Note that the upper/lower signs in \eqref{Kernel-K1plus} are the result of choosing to evaluate the factors at $z=1$ and $-1$ respectively. Ideally, the optimal choice for these constants would require the quantity $\widetilde{K}^+(z)-\widetilde{K}^-(1/z)$ to be minimised for all $z$ on the unit circle. However, experience with test cases have shown that optimisation makes a negligible difference to the value of the constants.

\begin{figure}[ht]\centering
\includegraphics[width=0.45\textwidth]{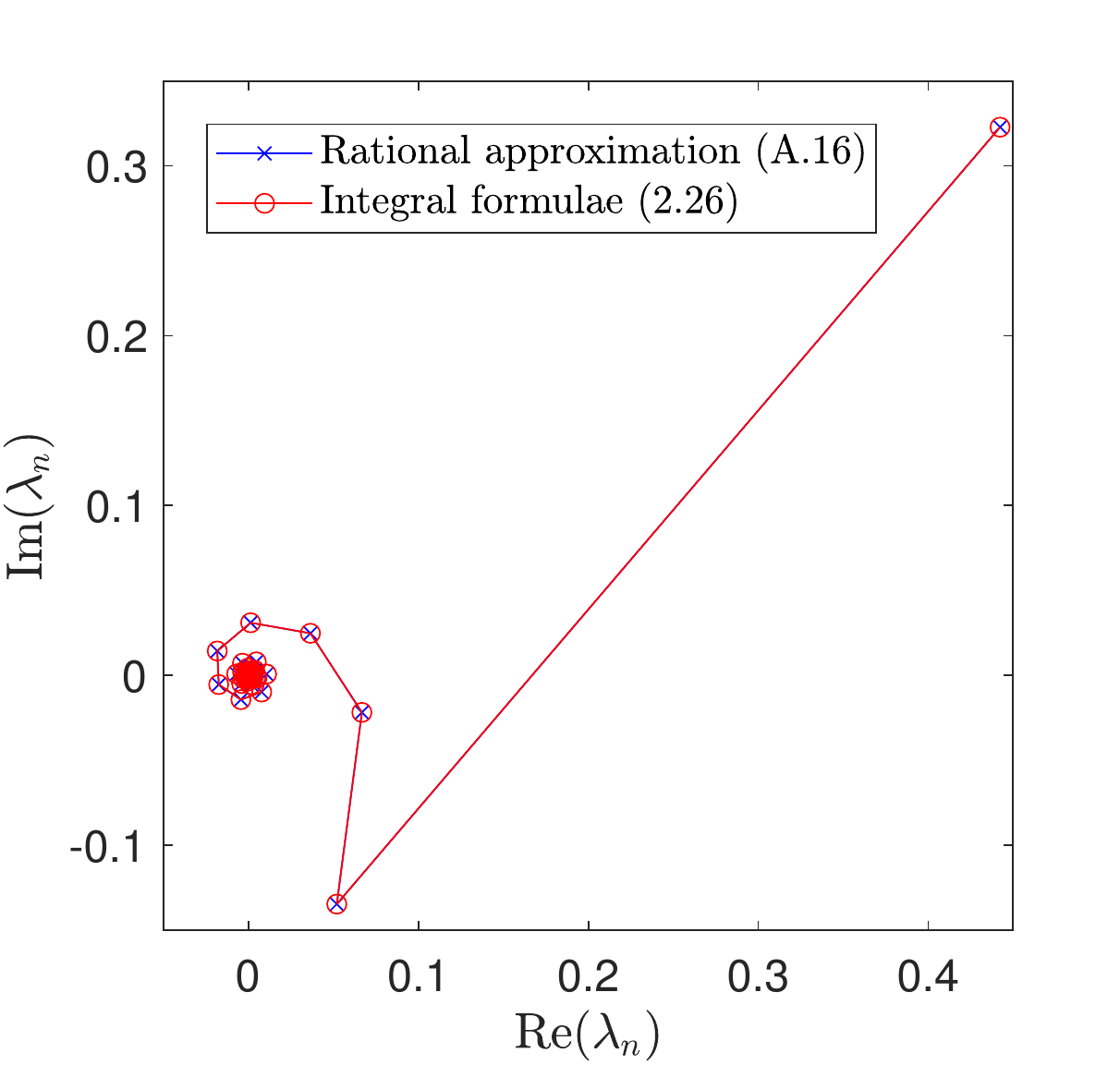}
\caption{A comparison of the $\lambda_n$ coefficients $(0\leq n\leq 1000)$ on the complex plane by the approximation formulae \eqref{Kernel-lambda-approx} ({\color{blue}blue}) and the integral formulae \eqref{SIG-lambda_n-integral} ({\color{red}red}). Note that the outlier is $\lambda_0$. Here we have $ks=1$ and $ka=0.01$.}
\label{fig:lambda_comparison}
\end{figure}

Lastly, we can use the rational approximation to accurately approximate the $\lambda_n$ coefficients. We perform a partial fraction expansion on the reciprocal of $\widetilde{K}^+(z)$,
\begin{align}\label{Kernel-K_approx-par-frac}
\frac{1}{\widetilde{K}^+(z)}&=\frac{1}{\widetilde{K}^+_1}\left[1+\sum_{\ell=1}^{L/2}\frac{\textrm{r}_\ell}{z-z^+_\ell}\right],\  
\ \textrm{where}\ \ \textrm{r}_m=(z^+_m-p^+_m)\prod_{\substack{\ell=1\\\ell\neq m}}^{L/2}\frac{z^+_m-p^+_\ell}{z^+_m-z^+_\ell}, 
\end{align}
and then we substitute \eqref{Kernel-K_approx-par-frac} into \eqref{SIG-1/K-series} to get these approximations for $\lambda_n$
\begin{align}\label{Kernel-lambda-approx}
\lambda_0&\approx\frac{1}{\widetilde{K}^+_1}\prod_{\ell=1}^{L/2}\frac{p^+_\ell}{z^+_\ell},\ \ \textrm{and}\ \ %
\lambda_n\approx-\frac{1}{\widetilde{K}^+_1}\sum_{\ell=1}^{L/2}\frac{\textrm{r}_\ell}{(z^+_\ell)^{n+1}}.
\end{align}
To test the accuracy of the $\lambda_n$ approximations, Figure \ref{fig:lambda_comparison} plots the $\lambda_n$ coefficients from $n=0$ to $1000$ on the complex plane using both the rational approximation \eqref{Kernel-lambda-approx} and the integral methods \eqref{SIG-lambda_n-integral}. 
We find that the two graphs are on top of each other. This is clear that evaluating the $\lambda_n$ coefficients via \eqref{Kernel-lambda-approx} is very accurate and we also find that it is significantly faster computationally than using the integral formula. However, it will only be as accurate as the rational approximation itself and numerical error will eventually catch up as $n\rightarrow\infty$ due to the denominator's exponent. 

\end{appendix}
\end{document}